# An Efficient, Second Order Accurate, Universal Generalized Riemann Problem Solver Based on the HLLI Riemann Solver

By

Dinshaw S. Balsara[1] , Jiequan Li[2] and Gino I. Montecinos[3]

[1]**Physics and ACMS Departments, University of Notre Dame (dbalsara@nd.edu)**

[2]**Institute of Applied Physics and Computational Mathematics, Beijing (li_jiequan@iapcm.ac.cn)**

[3]**University of Chile, Santiago (gmontecinos@cmm.dim.cl)**

**Abstract**

The Riemann problem, and the associated generalized Riemann problem, are increasingly seen as the important building blocks for modern higher order Godunov-type schemes. In the past, building a generalized Riemann problem solver was seen as an intricately mathematical task because the associated Riemann problem is different for each hyperbolic system of interest. This paper changes that situation.

The HLLI Riemann solver is a recently-proposed Riemann solver that is universal in that it is applicable to any hyperbolic system, whether in conservation form or with non-conservative products. The HLLI Riemann solver is also complete in the sense that if it is given a complete set of eigenvectors, it represents all waves with minimal dissipation. It is, therefore, very attractive to build a generalized Riemann problem solver version of the HLLI Riemann solver. This is the task that is accomplished in the present paper. We show that at second order, the generalized Riemann problem version of the HLLI Riemann solver is easy to design. Our GRP solver is also complete and universal because it inherits those good properties from original HLLI Riemann solver. We also show how our GRP solver can be adapted to the solution of hyperbolic systems with stiff source terms.

Our generalized HLLI Riemann solver is easy to implement and performs robustly and well over a range of test problems. All implementation-related details are presented. Results from several stringent test problems are shown. These test problems are drawn from many different hyperbolic systems, and include hyperbolic systems in conservation form; with non-conservative products; and with stiff source terms. The present generalized Riemann problem solver performs well on all of them.

# I) Introduction

The Riemann solver has been a building block of Godunov-type schemes since their inception (Godunov [31]). In the Riemann problem, two constant states coming together at a zone boundary provide the input variables for the Riemann solver. The resolved state and flux from the Riemann solver are the desired output variables that we seek from the Riemann solver. Kolgan [36] and van Leer [54] presented the first functional second order Godunov-type scheme. It was based on endowing each zone with a linear profile. Consequently, the input variables for the Riemann solver were not just states from both sides of a zone boundary but also their gradients in one-dimension. The Riemann solver then has to evaluate not just the resolved state and flux but also their variation with time in order to produce a time-centered, second order flux at the given zone boundary. This explains to us the concept of a generalized Riemann solver. The generalized Riemann problem solver can be characterized as a machine that accepts left and right states along with their higher order gradients at a zone boundary. It produces as output the resolved state and the resolved flux and also their time-evolution to the desired order of accuracy.

Progress in the design of second order accurate GRP solvers picked up with the work of Ben-Artzi and Falcovitz [8], [9] and Ben-Artzi [10]. A generalization to any general hyperbolic system in self-similar variables was made by LeFloch and Raviart [6] and Bourgeade *et al*. [18]. However, the work of LeFloch and Raviart is more of a systematization of GRP methodology. The textbook by Ben-Artzi and Falcovitz [12] also helped popularize GRP solvers. The field of GRP design got its second spurt of momentum from the work of Ben-Artzi, Li and Warnecke [13] for the Euler system and Ben-Artzi and Li [14] for general systems. A GRP for shallow water equations was also constructed by Li and Chen [38]. In Qian, Li and Wang [45] a third order accurate GRP solver was also designed for general hyperbolic systems and in Wu, Yang and Tang [56] for the Euler system. For the exact GRP for the Euler system, it has proved to be quite a daunting proposition to go past third order of accuracy. The application of the GRP solver to relativistic fluid dynamics was due to Yang and Tang [57], Wu and Tang [55]. As shown in Li and Wang [41] the analytic and nonlinear derivation of the GRP solver can cope with very extreme challenges such as high temperature and large ratio of density for multi-material flows because thermodynamical effects are precisely merged into the design of the solver.

Approximate solutions to the GRP were also attempted under the rubric of ADER schemes by Toro *et al*. [53], Titarev and Toro [51], [52], Montecinos and Toro [44]. An ADER scheme was also constructed for the MHD system by Taube et al. [50]. Goetz and Iske [32] and Goetz and Dumbser [33] also designed a GRP that was based on the method of LeFloch and Raviart. The second order ADER is the acoustic version of GRP in Ben-Artzi and Falcovitz [12], Ben-Artzi, Li and Warnecke [13], Ben-Artzi and Li [14] and Han *et al*. [35]. It is fair to observe that with a few exceptions most of these GRP solvers have been restricted to the Euler system. This is because the analytic constructions for the GRP in one way or the other involve the Cauchy-Kovalevskaya procedure. Alternatively, they require a detailed solution of the variation of the flow variables within a rarefaction fan or at a shock. The analytical skills that are involved in the construction of the previously-mentioned GRP solvers have always been a natural barrier in this sort of work.

Many of the above works attempted to construct exact GRP solvers for the associated Riemann problems. This proves to be a formidable task, especially when the hyperbolic system becomes very large, although they can be done in principle (see Ben-Artzi and Li [14], Qian, Li and Wang [45] Li and Sun [40]). One can always ask whether there is a general-purpose strategy for constructing approximate GRPs for entire classes of hyperbolic systems? Based on the application of approximate shock jumps at the wave boundaries of approximate Riemann solvers, Goetz, Balsara and Dumbser ([34]; GBD hereafter) showed that it is possible to design an approximate GRP solver out of the HLL, HLLC and HLLD Riemann solvers. The spatial gradients in the resolved states were found by a least squares procedure. This provided a considerable simplification in the construction of GRP solvers. The HLL-GRP solver is applicable to all manner of conservation laws, but it is rather dissipative in its treatment of intermediate waves. The HLLC-GRP only applies to the Euler system. The HLLD-GRP only applies to the MHD system. In summary, while the work of Goetz, Balsara and Dumbser [34] made it possible to apply the GRP philosophy to several approximate Riemann solvers, it was still specific to a small set of hyperbolic systems. Obtaining a general-purpose GRP solver that applies to any hyperbolic system, whether in conservation form or in non-conservative form, still eluded Goetz, Balsara and Dumbser [34]. This was especially true if one wished to resolve intermediate waves in the hyperbolic system.

What is truly desired is to start with an approximate Riemann solver that is universal and complete and build a GRP version of it. By universal we mean that it applies to any hyperbolic system, whether it is in conservation form or it has non-conservative products. By complete we mean that it resolves the full family of intermediate waves. An approximate Riemann solver that is universal and complete has been presented by Dumbser and Balsara [28]. This Riemann solver is called HLLI because it is built on top of the HLL Riemann solver and it can handle any number of intermediate waves; hence the name "HLLI". The HLLI Riemann solver follows in the path of the HLLEM Riemann solver (Einfeldt [29], Einfeldt *et al*. [30]) and was also motivated by recent self-similar formulations of the multidimensional Riemann solver (Balsara [4], Balsara and Dumbser [5], Balsara and Nkonga [6]). The *first goal* of this paper is to present a second order accurate GRP version of the HLLI Riemann solver of Dumbser and Balsara [28]; we refer to this as the HLLI-GRP solver.

The strategy used in the design of the HLLI-GRP in this paper builds on the approach of GBD. Thus, given states and their higher derivatives on either side of a zone boundary, GBD were able to obtain derivatives of the resolved HLL state within the Riemann fan. This was accomplished with a least squares procedure for conservation laws. In this paper our *second goal* is to extend that concept to include hyperbolic systems with non-conservative products. Thus a generalized form of shock jump conditions are derived and those jump conditions are used to obtain the gradient of the resolved state within the Riemann fan even when the hyperbolic system has non-conservative products. We then extend the HLLI Riemann solver construction so as to make it second order in time. We, therefore, obtain a flux form HLLI-GRP that applies to conservative systems. However, we also obtain a fluctuation form HLLI-GRP that applies to hyperbolic systems with non-conservative products.

The importance of the GRP in tackling problems involving stiff source terms has also been discussed by BenArtzi [10] and BenArtzi and Birman [11]. In recent years, this importance has also been stressed by Montecinos and Toro [44] and Goetz and Dumbser [33]. Such stiff sources arise in reactive flow calculations where the right temperature can initiate combustion in a large reaction network. Such terms also occur in radiation hydrodynamics, which seeks to couple the radiation field with the hydrodynamic equations. In all such situations, one seeks a GRP solver that is stiffly stable. In this paper our *third goal* is to address the treatment of stiff source terms in the HLL-GRP and HLLI-GRP solvers. Our approach is based on an innovative ADER scheme that is reported in Balsara *et al*. [7]. We show how this ADER scheme can be adapted so that it can be invoked *within* the Riemann fan! The result is that the HLL-GRP and HLLI-GRP solvers inherit all the good A-stable properties of the ADER scheme that is embedded in the GRP solution process.

The present paper on the HLL-GRP and HLLI-GRP solvers is generally useful for any scheme that wants to reduce the number of stages in the Runge-Kutta procedure. Such a step overcomes many limitations of the Butcher barriers that have plagued Strong Stability Preserving Runge-Kutta (SSP-RK) schemes. Consequently, Li and Du [39] presented a two stage scheme that is fourth order accurate in space and time. This is done by relying on a second order GRP. Likewise, Christlieb *et al*. [23] presented analogous two stage schemes that are potentially fourth order accurate. These schemes, therefore, occupy an intermediate position between the SSP-RK schemes (Shu and Osher [48], [49]) and the modern ADER schemes (Dumbser *et al*. [25], Balsara *et al*. [1], [3]). They are easier to implement than ADER schemes, yet they do not need to repeat the reconstruction step and Riemann solvers at so many sub-stages like the SSP-RK scheme. GRP solvers have also been used to guide the motion of the mesh in ALE codes (Boscheri, Balsara and Dumbser [16], Boscheri Dumbser and Balsara [17]). It is not the purpose of this paper to document the many advanced applications of the GRP that we design here. The goal of this paper is to thoroughly document the HLL-GRP and HLLI-GRP solvers for conservative and non-conservative systems and to show that they work well on a range of hyperbolic systems.

The plan of this paper is as follows. In Section II we describe the HLLI-GRP solver for conservation laws at second order. In Section III we formulate the HLL-GRP and HLLI-GRP solvers for hyperbolic systems with non-conservative products. We also describe the process for obtaining the gradient of the resolved state in the Riemann solver when non-conservative products are present. In Section IV we describe the HLLI-GRP-based scheme for hyperbolic systems in conservation form as well as hyperbolic systems with non-conservative products. In Section V we show how stiff source terms can be included in our formalism. Section VI presents several test problems. Section VII draws some conclusions.

## II) An HLLI-GRP Solver for Conservation Laws at Second Order

In this section we focus on a derivation of the HLLI-GRP solver for conservation laws at second order of accuracy. The derivation is split into three sub-sections. In Sub-section II.a we cast the conservation law into similarity variables and find a solution within each constant state. In Sub-section II.b we show how this contributes to the formulation of an HLL-GRP for conservation laws. In Sub-section II.c we show how this can be transcribed to the formulation of an HLLI-GRP for conservation laws at second order.

### II.a) One-Dimensional Conservation Law in Similarity Variables

Consider the $M \times M$ system of conservation laws in one-dimension of the form

$$\frac{\partial \mathbf{U}}{\partial t} + \frac{\partial \mathbf{F}(\mathbf{U})}{\partial x} = 0 \tag{2.1}$$

Here "$\mathbf{U}$" is an $M$-component vector and "$\mathbf{F}(\mathbf{U})$" is an $M$-component flux that depends on "$\mathbf{U}$". Since the solution to the Riemann problem is self-similar, it pays to cast the above equation in similarity variables. Following LeFloch and Raviart [37], as well as Ben-Artzi and Li [14] for the acoustic version, we take the similarity variable $\xi \equiv x/t$ and make the change of variables $(x,t) \to (\xi,t)$ to recast eqn. (2.1) in the form

$$t\frac{\partial \mathbf{U}}{\partial t} - \xi\frac{\partial \mathbf{U}}{\partial \xi} + \frac{\partial \mathbf{F}(\mathbf{U})}{\partial \xi} = 0 \tag{2.2}$$

The GRP basically incorporates the deviations from self-similarity in terms of a power series in time around $t = 0$. We assert a Taylor series in time solution of the form

$$\mathbf{U}(\xi,t) = \mathbf{U}^0(\xi) + t\,\mathbf{U}^1(\xi) + \frac{1}{2}t^2\,\mathbf{U}^2(\xi) + \ldots + \frac{1}{k!}t^k\,\mathbf{U}^k(\xi) + \ldots \tag{2.3}$$

Since it is our goal to start by finding a solution to the HLL-GRP, which starts with a constant resolved state, we take $\mathbf{U}^0(\xi) \equiv \mathbf{U}^0$ to be a constant and only carry the first two terms in eqn. (2.3).

Defining $\mathbf{A}(\mathbf{U}^0) = \partial \mathbf{F}(\mathbf{U}^0) / \partial \mathbf{U}$, we obtain an equation for the flux and its linearized variations with respect to variations in the state as

$$\mathbf{F}(\xi,t) = \mathbf{F}(\mathbf{U}^0) + \mathbf{A}(\mathbf{U}^0)(\mathbf{U}(\xi,t) - \mathbf{U}^0) \qquad (2.4)$$

Using eqns. (2.3) and (2.4), eqn. (2.2) can be specialized to yield an equation for $\mathbf{U}^1(\xi)$. The resulting equation is

$$\mathbf{U}^1 - (\xi I - \mathbf{A}(\mathbf{U}^0)) \frac{\partial \mathbf{U}^1}{\partial \xi} = 0 \qquad (2.5)$$

We see, therefore, that $\mathbf{U}^1(\xi)$ will be linear in the variable "$\xi$" and, therefore, the solutions for eqns. (2.3) and (2.4), when retaining terms that are linear in time "$t$", are given by

$$\tilde{\mathbf{U}}(\xi,t) = \mathbf{U}^0(\xi) + t\ (\xi \mathbf{I} - \mathbf{A}(\mathbf{U}^0))(\partial_x \mathbf{U}^1) \qquad (2.6)$$

and

$$\tilde{\mathbf{F}}(\xi,t) = \mathbf{F}(\mathbf{U}^0) + t\ \mathbf{A}(\mathbf{U}^0)(\xi \mathbf{I} - \mathbf{A}(\mathbf{U}^0))(\partial_x \mathbf{U}^1) \qquad (2.7)$$

Here $(\partial_x \mathbf{U}^1)$ is an *M*-component vector that is identified with the first derivative of the solution in space. Eqns. (2.6) and (2.7) give us the two essential equations that we will use in the next section. The tilde on top of $\tilde{\mathbf{U}}(\xi,t)$ and $\tilde{\mathbf{F}}(\xi,t)$ is intended to highlight the functional nature of the state and flux, and this is a notation that we will use consistently in this paper.

### II.b) Formulation of the HLL-GRP Solver for Conservation Laws at Second Order

We start by describing the input states to the GRP. We have the left state $\mathbf{U}^L$ and its spatial derivative $(\partial_x \mathbf{U}^L)$ and the right state $\mathbf{U}^R$ and its spatial derivative $(\partial_x \mathbf{U}^R)$. The states $\mathbf{U}^L$ and $\mathbf{U}^R$ can be used to identify an extremal left-going wave with speed $S_L$ and an extremal right-going wave with speed $S_R$. In a space-time diagram, the HLL Riemann solver can be obtained by asserting that these extremal waves contain a resolved state $\mathbf{U}^*$ and a resolved flux $\mathbf{F}^*$, which are given by

$$\mathbf{U}^* = \frac{S_R \mathbf{U}^R - S_L \mathbf{U}^L - (\mathbf{F}(\mathbf{U}^R) - \mathbf{F}(\mathbf{U}^L))}{(S_R - S_L)} \qquad (2.8)$$

and

$$\mathbf{F}^* = \frac{S_R \mathbf{F}(\mathbf{U}^L) - S_L \mathbf{F}(\mathbf{U}^R) + S_R S_L (\mathbf{U}^R - \mathbf{U}^L)}{(S_R - S_L)} \tag{2.9}$$

Our goal in studying the HLL-GRP is to endow the resolved state $\mathbf{U}^*$ with a spatial gradient $(\partial_x \mathbf{U}^*)$. Once that is done, as will be shown in the ensuing paragraph, we may provide the space-time evolution of the resolved state for $S_L < \xi < S_R$ in the $t \to 0$ limit as follows

$$\tilde{\mathbf{U}}^*(\xi, t) = \mathbf{U}^* + t \left( \xi \mathbf{I} - \mathbf{A}(\mathbf{U}^*) \right) (\partial_x \mathbf{U}^*) \tag{2.10}$$

and

$$\tilde{\mathbf{F}}^*(\xi, t) = \mathbf{F}^* + t \, \mathbf{A}(\mathbf{U}^*) \left( \xi \mathbf{I} - \mathbf{A}(\mathbf{U}^*) \right) (\partial_x \mathbf{U}^*) \tag{2.11}$$

This is done by imposing shock jump conditions at both the extremal waves, as shown in the next paragraph.

We illustrate this process by considering the shock jump conditions at the right extremal wave. Analogous to eqns. (2.10) and (2.11), we have for $\xi \geq S_R$ in the $t \to 0$ limit

$$\tilde{\mathbf{U}}^R(\xi, t) = \mathbf{U}^R + t \left( \xi \mathbf{I} - \mathbf{A}(\mathbf{U}^R) \right) (\partial_x \mathbf{U}^R) \tag{2.12}$$

and

$$\tilde{\mathbf{F}}^R(\xi, t) = \mathbf{F}(\mathbf{U}^R) + t \, \mathbf{A}(\mathbf{U}^R) \left( \xi \mathbf{I} - \mathbf{A}(\mathbf{U}^R) \right) (\partial_x \mathbf{U}^R) \tag{2.13}$$

We also make allowance for the fact that the extremal right-going wave may have a curved trajectory in space-time that is given by

$$x_R(t) = S_R t + \frac{1}{2} S_R' t^2 \tag{2.14}$$

In the above equation, $S_R'$ carries the curvature of the wave in space-time. The Rankine-Hugoniot shock-jump condition, extended to the GRP, and asserted at the right-going extremal wave of the HLL Riemann solver then reads

$$\tilde{\mathbf{F}}^*(\xi = S_R, t) - \tilde{\mathbf{F}}^R(\xi = S_R, t) = \frac{dx_R(t)}{dt} \left[ \tilde{\mathbf{U}}^*(\xi = S_R, t) - \tilde{\mathbf{U}}^R(\xi = S_R, t) \right] \tag{2.15}$$

By focusing on the time-dependent part of the above right-going shock-jump condition we get

$$\left( S_R \mathbf{I} - \mathbf{A}(\mathbf{U}^*) \right)^2 (\partial_x \mathbf{U}^*) + S_R' (\mathbf{U}^* - \mathbf{U}^R) = \left( S_R \mathbf{I} - \mathbf{A}(\mathbf{U}^R) \right)^2 (\partial_x \mathbf{U}^R) \tag{2.16}$$

An analogous condition at the left-going extremal wave of the HLL Riemann solver then gives us

$$\left(S_L\mathbf{I}-\mathbf{A}\left(\mathbf{U}^*\right)\right)^2\left(\partial_x\mathbf{U}^*\right)+S_L'\left(\mathbf{U}^*-\mathbf{U}^L\right)=\left(S_L\mathbf{I}-\mathbf{A}\left(\mathbf{U}^L\right)\right)^2\left(\partial_x\mathbf{U}^L\right) \quad (2.17)$$

Eqns. (2.16) and (2.17) are extremely important for the GRP because they give us $\left(\partial_x\mathbf{U}^*\right)$ which leads us to the solution of the GRP for conservative hyperbolic systems. The unknowns in the above two equations consist of the *M*-component vector $\left(\partial_x\mathbf{U}^*\right)$ and the two scalars $S_R'$ and $S_L'$. Consequently, we have $(M+2)$ unknowns whereas eqns. (2.16) and (2.17) constitute $2M$ linear equations in those unknowns. Therefore, for $M\geq 2$, a solution can always be found by minimizing the variation in eqns. (2.16) and (2.17) in a least squares sense. This completes the description our solution methodology for obtaining the *M*-component gradient vector in the resolved state, $\left(\partial_x\mathbf{U}^*\right)$, as well as the curvature terms $S_R'$ and $S_L'$. By obtaining the above-mentioned $(M+2)$ unknowns we solve the HLL-GRP for conservation laws.

If one does not care to retain the curvature terms in the extremal waves, eqns. (2.16) and (2.17) can also be viewed as $2M$ equations in $M$ unknowns. This was the choice made in GBD. We have experimented with that possibility and found it to work very well too. For that reason, we used this approximation in all our numerical experiments.

### II.c) Extension to the HLLI-GRP Solver for Conservation Laws at Second Order

Realize that the solution of eqns. (2.16) and (2.17) gives us $\left(\partial_x\mathbf{U}^*\right)$ which is the gradient of the resolved state in the HLL Riemann solver. Transitioning from the HLL-GRP to the HLLI-GRP is indeed quite simple at second order. At second order, there are no cross terms between space and time. Say we wish to take a timestep of size $\Delta t$ from time $t^n$ to time $t^{n+1}$ so that we want to evaluate time-centered fluxes at $\Delta t/2$. We can then evaluate the contribution from the anti-diffusive fluxes at time $\Delta t/2$ by evaluating the right and left states at an advanced time. Denote them by $\mathbf{U}^{R;n+1/2}$ and $\mathbf{U}^{L;n+1/2}$ respectively. They are obtained as follows

$$\mathbf{U}^{R;n+1/2}\equiv\tilde{\mathbf{U}}^R\left(\xi=0,t=\Delta t/2\right)=\mathbf{U}^R-\frac{\Delta t}{2}\mathbf{A}\left(\mathbf{U}^R\right)\left(\partial_x\mathbf{U}^R\right) \quad (2.18)$$

and

$$\mathbf{U}^{L;n+1/2}\equiv\tilde{\mathbf{U}}^L\left(\xi=0,t=\Delta t/2\right)=\mathbf{U}^L-\frac{\Delta t}{2}\mathbf{A}\left(\mathbf{U}^L\right)\left(\partial_x\mathbf{U}^L\right) \quad (2.19)$$

We can also obtain our HLL state and flux at an advanced time. Denote them by $\mathbf{U}^{*;n+1/2}$ and $\mathbf{F}_{HLL-GRP}$ respectively. They are obtained as follows

$$\mathbf{U}^{*;n+1/2}\equiv\tilde{\mathbf{U}}^*\left(\xi=0,t=\Delta t/2\right)=\mathbf{U}^*-\frac{\Delta t}{2}\mathbf{A}\left(\mathbf{U}^*\right)\left(\partial_x\mathbf{U}^*\right) \quad (2.20)$$

and

$$\mathbf{F}_{HLL-GRP} \equiv \tilde{\mathbf{F}}^*\left(\xi = 0, t = \Delta t/2\right) = \mathbf{F}^* - \frac{\Delta t}{2}\left(\mathbf{A}\left(\mathbf{U}^*\right)\right)^2\left(\partial_x \mathbf{U}^*\right) \quad (2.21)$$

To obtain the HLLI-GRP we have to first decide on the set of waves that we want to improve. For an $M \times M$ hyperbolic system, we identify "$P$" waves that we want to improve, where $P \leq M$. Using $\mathbf{U}^{*;n+1/2}$, or by using $\left(\mathbf{U}^{R;n+1/2} + \mathbf{U}^{L;n+1/2}\right)/2$, we can evaluate the set of eigenvalues $\{\lambda^p; p = 1,..,P\}$, the corresponding set of left eigenvectors $\{l^p; p = 1,..,P\}$ and the corresponding set of right eigenvectors $\{r^p; p = 1,..,P\}$. The HLLI-GRP flux can now be written in terms of the HLL-GRP flux as follows

$$\mathbf{F}_{HLLI-GRP} = \mathbf{F}_{HLL-GRP} - \phi \frac{S_R S_L}{(S_R - S_L)} \sum_{p=1}^{P} \delta^p \left[l^p \cdot \left(\mathbf{U}^{R;n+1/2} - \mathbf{U}^{L;n+1/2}\right)\right] r^p \quad (2.22)$$

Here $\phi$ is a flattener function (Colella and Woodward [19], Balsara [2]) which is unity away from shocks and becomes closer to zero at strong shocks. It can smoothly assume any value between zero and one, so that shocks of intermediate strength can be properly treated. In other words, we are following the suggestion of Dumbser and Balsara [28] who find that the anti-diffusive contribution should be suppressed at shocks. Following Dumbser and Balsara [28] we set

$$\delta^p = 1 - \frac{\min(\lambda^p, 0)}{S_L} - \frac{\max(\lambda^p, 0)}{S_R} \quad (2.23)$$

This completes our description of the HLLI-GRP for conservative systems. The insights developed here will be very valuable for understanding the GRP for non-conservative systems, which we study next.

### III) An HLLI-GRP Solver for Hyperbolic Systems with Non-Conservative Products at Second Order

In this section we focus on a derivation of the HLLI-GRP for hyperbolic systems with non-conservative products at second order of accuracy. The derivation is split into three sub-sections. In Sub-section III.a we cast the hyperbolic system into similarity variables and find a solution within each constant state. In Sub-section III.b we show how this contributes to the formulation of an HLL Riemann solver for hyperbolic systems with non-conservative products. In Sub-section III.c we show how to obtain the gradient terms in the resolved state for the HLL-GRP. The complete assembly of an HLL-GRP or an HLLI-GRP will be presented in Section IV.

### III.a) One-Dimensional Hyperbolic System with Non-Conservative Products in Similarity Variables

Consider the one-dimensional $M \times M$ hyperbolic system with non-conservative products as follows

$$\frac{\partial \mathbf{U}}{\partial t} + \frac{\partial \mathbf{F}(\mathbf{U})}{\partial x} + \mathbf{B}(\mathbf{U})\frac{\partial \mathbf{U}}{\partial x} = 0 \qquad (3.1)$$

Here the $M \times M$ matrix "$\mathbf{B}(\mathbf{U})$" contains the non-conservative products. As before, "$\mathbf{U}$" is an $M$-component vector and "$\mathbf{F}(\mathbf{U})$" is an $M$-component flux that depends on "$\mathbf{U}$". As before, we cast the system in terms of similarity variables so that the transcription $(x,t) \rightarrow (\xi,t)$ gives us

$$t\frac{\partial \mathbf{U}}{\partial t} + \mathbf{U} - \frac{\partial(\xi \mathbf{U})}{\partial \xi} + \frac{\partial \mathbf{F}(\mathbf{U})}{\partial \xi} + \mathbf{B}(\mathbf{U})\frac{\partial \mathbf{U}}{\partial \xi} = 0 \qquad (3.2)$$

We assert a Taylor series in time solution of the form given by eqn. (2.3). Defining $\mathbf{C}(\mathbf{U}) = \partial \mathbf{F}(\mathbf{U})/\partial \mathbf{U}$, the characteristic matrix for the system can be written as $\mathbf{A}(\mathbf{U}) = \mathbf{C}(\mathbf{U}) + \mathbf{B}(\mathbf{U})$. Eqn. (2.5) is, therefore, unchanged with the result that eqn. (2.6) for $\tilde{\mathbf{U}}(\xi,t)$ is still valid. Eqn. (2.7) which gives $\tilde{\mathbf{F}}(\xi,t)$ for the flux does, however, change so that we have

$$\tilde{\mathbf{F}}(\xi,t) = \mathbf{F}(\mathbf{U}^0) + t\, \mathbf{C}(\mathbf{U}^0)\bigl(\xi \mathbf{I} - \mathbf{A}(\mathbf{U}^0)\bigr)(\partial_x \mathbf{U}^1) \qquad (3.3)$$

Here $(\partial_x \mathbf{U}^1)$ is an $M$-component vector that is identified with the first derivative of the solution in space. Eqns. (2.6) and (3.3) give us the two essential equations that we will use in the next section.

### III.b) Formulation of the HLL Riemann Solver for Hyperbolic Systems with Non-Conservative Products

In Dumbser and Balsara [28] we provided a derivation of the HLL Riemann solver for non-conservative hyperbolic systems that was based on Heaviside weights. The use of Heaviside weights may not be transparent for all readers. Let us, therefore, provide a much-simplified derivation of the resolved HLL state $\mathbf{U}^*$ for the HLL Riemann problem when non-conservative products are present. The extremal wave speeds $S_L$ and $S_R$ are obtained as in the previous section. The states in the HLL Riemann solver can be parametrized along the following path

$$\overline{\mathbf{U}}(\xi) = \begin{cases} \mathbf{U}^L & \xi \leq S_L - \varepsilon \\ \mathbf{U}^L + (\mathbf{U}^* - \mathbf{U}^L)\dfrac{(\xi - S_L + \varepsilon)}{\varepsilon} & S_L - \varepsilon < \xi < S_L \\ \mathbf{U}^* & S_L \leq \xi \leq S_R \\ \mathbf{U}^* + (\mathbf{U}^R - \mathbf{U}^*)\dfrac{(\xi - S_R)}{\varepsilon} & S_L - \varepsilon < \xi < S_L \\ \mathbf{U}^R & S_R + \varepsilon \leq \xi \end{cases} \quad (3.4)$$

In eqn. (3.4) we parametrize the path $\mathbf{U}^L \to \mathbf{U}^* \to \mathbf{U}^R$ directly in terms of the similarity variable as shown in Fig. 1. In the limit $\varepsilon \to 0$, eqn. (3.4) gives an equivalent derivation of the HLL state $\mathbf{U}^*$ to the derivation in Dumbser and Balsara [28]. However, the present derivation may be more intuitively appealing to some people because it does not rely on the use of Heaviside weights. Since the HLL Riemann solver (i.e. without the GRP part) has no time evolution, we can take the time-independent version of eqn. (3.2) and integrate it from $\xi = S_L - \varepsilon$ to $\xi = S_R + \varepsilon$. We, therefore, evaluate

$$\int_{\xi = S_L - \varepsilon}^{\xi = S_R + \varepsilon} \left( \overline{\mathbf{U}}(\xi) - \frac{\partial(\xi \overline{\mathbf{U}}(\xi))}{\partial \xi} + \frac{\partial \mathbf{F}(\overline{\mathbf{U}}(\xi))}{\partial \xi} + \mathbf{B}(\overline{\mathbf{U}}(\xi))\frac{\partial \overline{\mathbf{U}}(\xi)}{\partial \xi} \right) d\xi = 0 \quad (3.5)$$

The specific form of $\overline{\mathbf{U}}(\xi)$ from eqn. (3.4) can be used in the above integral. The result, after neglecting terms that are of order "$\varepsilon$", is

$$\mathbf{U}^*(S_R - S_L) - (S_R \mathbf{U}^R - S_L \mathbf{U}^L) + (\mathbf{F}(\mathbf{U}^R) - \mathbf{F}(\mathbf{U}^L))$$
$$+ \int_{\xi = S_L - \varepsilon}^{\xi = S_L} \mathbf{B}(\overline{\mathbf{U}}(\xi))\frac{\partial \overline{\mathbf{U}}(\xi)}{\partial \xi} d\xi + \int_{\xi = S_R}^{\xi = S_R + \varepsilon} \mathbf{B}(\overline{\mathbf{U}}(\xi))\frac{\partial \overline{\mathbf{U}}(\xi)}{\partial \xi} d\xi = 0 \quad (3.6)$$

The two integrals in the above equation can be interpreted like a parameter vector style evaluation of the Roe matrix. Thus we can define the matrices $\tilde{\mathbf{B}}(\mathbf{U}^L, \mathbf{U}^*)$ and $\tilde{\mathbf{B}}(\mathbf{U}^*, \mathbf{U}^R)$ via the following two equations which use the parametrization in eqn. (3.4) to give

$$\int_{\xi = S_L - \varepsilon}^{\xi = S_L} \mathbf{B}(\overline{\mathbf{U}}(\xi))\frac{\partial \overline{\mathbf{U}}(\xi)}{\partial \xi} d\xi = \left[ \int_{\xi = S_L - \varepsilon}^{\xi = S_L} \mathbf{B}(\overline{\mathbf{U}}(\xi)) d\xi \right](\mathbf{U}^* - \mathbf{U}^L) \equiv \tilde{\mathbf{B}}(\mathbf{U}^L, \mathbf{U}^*)(\mathbf{U}^* - \mathbf{U}^L) \quad (3.7)$$

and

$$\int_{\xi = S_R}^{\xi = S_R + \varepsilon} \mathbf{B}(\overline{\mathbf{U}}(\xi))\frac{\partial \overline{\mathbf{U}}(\xi)}{\partial \xi} d\xi = \left[ \int_{\xi = S_R}^{\xi = S_R + \varepsilon} \mathbf{B}(\overline{\mathbf{U}}(\xi)) d\xi \right](\mathbf{U}^R - \mathbf{U}^*) \equiv \tilde{\mathbf{B}}(\mathbf{U}^*, \mathbf{U}^R)(\mathbf{U}^R - \mathbf{U}^*) \quad (3.8)$$

The integrals in the square brackets in eqns. (3.7) and (3.8) serve to define the matrices $\tilde{\mathbf{B}}(\mathbf{U}^L, \mathbf{U}^*)$ and $\tilde{\mathbf{B}}(\mathbf{U}^*, \mathbf{U}^R)$ respectively. In practice, the integrals in the square brackets in the above two equations are evaluated via a sufficiently accurate numerical quadrature. In other words, we evaluate the integrals in exactly the same fashion as for the DOT Riemann solver ( Dumbser and Toro [26], [27]). Eqn. (3.6) can now be written as

$$\mathbf{U}^* = \frac{1}{(S_R - S_L)} \left[ (S_R \mathbf{U}^R - S_L \mathbf{U}^L) - (\mathbf{F}(\mathbf{U}^R) - \mathbf{F}(\mathbf{U}^L)) - \tilde{\mathbf{B}}(\mathbf{U}^L, \mathbf{U}^*)(\mathbf{U}^* - \mathbf{U}^L) - \tilde{\mathbf{B}}(\mathbf{U}^*, \mathbf{U}^R)(\mathbf{U}^R - \mathbf{U}^*) \right]$$

(3.9)

The above equation is implicit in $\mathbf{U}^*$ and Dumbser and Balsara [28] provide easily-implemented methods for its efficient solution via an inexact Newton procedure. While this derivation of the state $\mathbf{U}^*$ may be easier to follow, it also serves another very useful purpose. In the next paragraph we will show that it leads us to a generalization of the shock-jump conditions for non-conservative systems.

For the HLL Riemann solver, the expressions for the left-going and right-going fluctuations are given by

$$\mathbf{D}_{HLL}^-(\mathbf{U}^L, \mathbf{U}^R) = \int_{\xi = S_L - \varepsilon}^{\xi = 0} \mathbf{A}(\bar{\mathbf{U}}(\xi)) \frac{\partial \bar{\mathbf{U}}(\xi)}{\partial \xi} d\xi = -\int_{\xi = S_L - \varepsilon}^{\xi = 0} \left( \bar{\mathbf{U}}(\xi) - \frac{\partial(\xi \bar{\mathbf{U}}(\xi))}{\partial \xi} \right) d\xi = S_L(\mathbf{U}^* - \mathbf{U}^L)$$

(3.10)

and

$$\mathbf{D}_{HLL}^+(\mathbf{U}^L, \mathbf{U}^R) = \int_{\xi = 0}^{\xi = S_R + \varepsilon} \mathbf{A}(\bar{\mathbf{U}}(\xi)) \frac{\partial \bar{\mathbf{U}}(\xi)}{\partial \xi} d\xi = -\int_{\xi = 0}^{\xi = S_R + \varepsilon} \left( \bar{\mathbf{U}}(\xi) - \frac{\partial(\xi \bar{\mathbf{U}}(\xi))}{\partial \xi} \right) d\xi = S_R(\mathbf{U}^R - \mathbf{U}^*)$$

(3.11)

The specific form of $\bar{\mathbf{U}}(\xi)$ from eqn. (3.4) can be used in the above integrals. In the previous two equations we also use the time-independent part of eqn. (3.2) to transition from the integrals over the characteristic matrix to the integrals over the states. The integrals over the states are much easier to evaluate, of course. It is easy to verify that for the specific path in phase space that is chosen in eqn. (3.4) it is easy to verify the consistency condition for the fluctuations

$$\mathbf{D}_{HLL}^-(\mathbf{U}^L, \mathbf{U}^R) + \mathbf{D}_{HLL}^+(\mathbf{U}^L, \mathbf{U}^R) = \int_{\xi = S_L - \varepsilon}^{\xi = S_R + \varepsilon} \mathbf{A}(\bar{\mathbf{U}}(\xi)) \frac{\partial \bar{\mathbf{U}}(\xi)}{\partial \xi} d\xi = \int_{\mathbf{U}^L}^{\mathbf{U}^R} \mathbf{A}(\bar{\mathbf{U}}) d\bar{\mathbf{U}}$$

(3.12)

This completes our description of the HLL Riemann solver for hyperbolic systems with non-conservative products.

For the sake of completeness, we also provide the formulae for the fluctuations when one goes from the HLL Riemann solver to the HLLI Riemann solver. The fluctuations are given by

$$\mathbf{D}_{HLLI}^{-}\left(\mathbf{U}^L,\mathbf{U}^R\right) = \mathbf{D}_{HLL}^{-}\left(\mathbf{U}^L,\mathbf{U}^R\right) - \phi\frac{S_R S_L}{(S_R - S_L)}\sum_{p=1}^{P}\delta^p\left[l^p\cdot\left(\mathbf{U}^R - \mathbf{U}^L\right)\right]r^p \quad (3.13)$$

and

$$\mathbf{D}_{HLLI}^{+}\left(\mathbf{U}^L,\mathbf{U}^R\right) = \mathbf{D}_{HLL}^{+}\left(\mathbf{U}^L,\mathbf{U}^R\right) + \phi\frac{S_R S_L}{(S_R - S_L)}\sum_{p=1}^{P}\delta^p\left[l^p\cdot\left(\mathbf{U}^R - \mathbf{U}^L\right)\right]r^p \quad (3.14)$$

As in the conservative case, we can evaluate the set of eigenvalues $\{\lambda^p; p=1,..,P\}$, the corresponding set of left eigenvectors $\{l^p; p=1,..,P\}$ and the corresponding set of right eigenvectors $\{r^p; p=1,..,P\}$ at some suitable average of $\mathbf{U}^L$ and $\mathbf{U}^R$. Notice that the anti-diffusive fluctuation terms make equal and opposite contributions to the fluctuations on either side of a zone boundary, as expected. This completes our description of the HLL Riemann solver for hyperbolic systems with non-conservative products.

### III.c) Obtaining the Gradient Terms in the Resolved State for the HLL-GRP Solver

We wish to obtain a GRP for non-conservative systems in a fashion that is analogous to the GRP in Sub-section II.b for conservative systems. For conservation laws, we have the benefit of having the Rankine-Hugoniot shock-jump conditions. We would like to generalize the shock-jump conditions so that, at least in a path conservative sense, they are also valid for a hyperbolic system with non-conservative products. Let us focus on the extremal right-going wave. In a distribution sense, it extends from $\xi = S_R$ to $\xi = S_R + \varepsilon$. The generalized shock-jump conditions can be obtained by integrating the time-independent version of eqn. (3.2) from $\xi = S_R$ to $\xi = S_R + \varepsilon$ as follows

$$\int_{\xi=S_R}^{\xi=S_R+\varepsilon}\left(\overline{\mathbf{U}}(\xi) - \frac{\partial(\xi\overline{\mathbf{U}}(\xi))}{\partial\xi} + \frac{\partial\mathbf{F}(\overline{\mathbf{U}}(\xi))}{\partial\xi} + \mathbf{B}(\overline{\mathbf{U}}(\xi))\frac{\partial\overline{\mathbf{U}}(\xi)}{\partial\xi}\right)d\xi = 0 \quad (3.15)$$

The specific form of $\overline{\mathbf{U}}(\xi)$ from eqn. (3.4) can be used in the above integral. Ignoring infinitesimal terms that vanish as $\varepsilon \to 0$ we get

$$\left(\mathbf{F}(\mathbf{U}^R) - S_R\mathbf{U}^R\right) - \left(\mathbf{F}^* - S_R\mathbf{U}^*\right) + \tilde{\mathbf{B}}\left(\mathbf{U}^*,\mathbf{U}^R\right)\left(\mathbf{U}^R - \mathbf{U}^*\right) = 0 \quad (3.16)$$

In the above equation, "$S_R$" refers to the instantaneous shock speed and it can be time-dependent in the GRP case. Similarly, the states in eqn. (3.16) can have time-dependence in the GRP case. Of course, we will rewrite the above equation so as to highlight this time-dependence for the GRP application in a short while; i.e., by the time we get to eqn. (3.20). The previous equation (along with the notational improvements that we introduce in eqn. (3.20) ) is indeed our master equation for shock jumps involving hyperbolic systems with non-conservative products. Using this equation, we will be able to assert our generalization of the Rankine-Hugoniot shock-jump condition. The above equation will enable us to relate gradients in the right state $\left(\partial_x \mathbf{U}^R\right)$ to gradients in the resolved state $\left(\partial_x \mathbf{U}^*\right)$ as will be shown in the next paragraph. For now, we also write a generalization of the Rankine-Hugoniot shock-jump condition for the extremal left-going wave. Asserting an equation that is analogous to eqn. (3.15) then gives us

$$\left(\mathbf{F}^* - S_L \mathbf{U}^*\right) - \left(\mathbf{F}\left(\mathbf{U}^L\right) - S_L \mathbf{U}^L\right) + \tilde{\mathbf{B}}\left(\mathbf{U}^L, \mathbf{U}^*\right)\left(\mathbf{U}^* - \mathbf{U}^L\right) = 0 \qquad (3.17)$$

In the above equation, "$S_L$" refers to the instantaneous shock speed and it can be time-dependent in the GRP case. The above equation will enable us to relate gradients in the left state $\left(\partial_x \mathbf{U}^L\right)$ to gradients in the resolved state $\left(\partial_x \mathbf{U}^*\right)$ as will be shown in the next paragraph. Armed with eqns. (3.16) and (3.17) we will derive the HLL-GRP for the case of hyperbolic systems with non-conservative products in the subsequent paragraphs.

Let us consider the right-going extremal wave. Eqn. (2.14) still gives us the trajectory of the right-going shock in a fashion that incorporates the curvature of the shock front. Eqns. (2.10) and (2.12) for the time-evolution of the states $\tilde{\mathbf{U}}^*(\xi,t)$ and $\tilde{\mathbf{U}}^R(\xi,t)$ respectively remain unchanged. In place of eqns. (2.11) and (2.13) for the evolution of the fluxes, we have

$$\tilde{\mathbf{F}}^*(\xi,t) = \mathbf{F}^* + t\ \mathbf{C}\left(\mathbf{U}^*\right)\left(\xi \mathbf{I} - \mathbf{A}\left(\mathbf{U}^*\right)\right)\left(\partial_x \mathbf{U}^*\right) \qquad (3.18)$$

and

$$\tilde{\mathbf{F}}^R(\xi,t) = \mathbf{F}\left(\mathbf{U}^R\right) + t\ \mathbf{C}\left(\mathbf{U}^R\right)\left(\xi \mathbf{I} - \mathbf{A}\left(\mathbf{U}^R\right)\right)\left(\partial_x \mathbf{U}^R\right) \qquad (3.19)$$

In place of the shock-jump condition in eqn. (2.15) we now use eqn. (3.16). When eqn. (3.16) is written in all generality, i.e. by making the time-dependence explicit for the GRP case, we get the following equation

$$\begin{aligned}&\tilde{\mathbf{F}}^R\left(\xi = S_R, t\right) - \tilde{\mathbf{F}}^*\left(\xi = S_R, t\right) - \frac{dx_R(t)}{dt}\left[\tilde{\mathbf{U}}^R\left(\xi = S_R, t\right) - \tilde{\mathbf{U}}^*\left(\xi = S_R, t\right)\right] \\ &+ \tilde{\mathbf{B}}\left(\tilde{\mathbf{U}}^*\left(\xi = S_R, t\right), \tilde{\mathbf{U}}^R\left(\xi = S_R, t\right)\right)\left(\tilde{\mathbf{U}}^R\left(\xi = S_R, t\right) - \tilde{\mathbf{U}}^*\left(\xi = S_R, t\right)\right) = 0\end{aligned} \qquad (3.20)$$

The above equation is the precise expression of the shock-jump condition for a hyperbolic system with non-conservative products when the GRP is being considered. The challenging part of the above equation consists of the last term involving the matrix $\tilde{\mathbf{B}}\left(\tilde{\mathbf{U}}^*\left(\xi = S_R, t\right), \tilde{\mathbf{U}}^R\left(\xi = S_R, t\right)\right)$. Notice that this Roe-linearization type of matrix depends on two states, $\tilde{\mathbf{U}}^*\left(\xi = S_R, t\right)$ and $\tilde{\mathbf{U}}^R(\xi = S_R, t)$ and will have to be expanded as such. The derivatives of this matrix, if taken seriously, will result in tensorial expressions (i.e., we will have to consider a Hessian). We carry these detailed terms for the sake of completeness. Thus we have the expansions up to *t*-dependent terms as

$$\tilde{\mathbf{B}}\left(\tilde{\mathbf{U}}^*(\xi = S_R, t), \tilde{\mathbf{U}}^R(\xi = S_R, t)\right)\left(\tilde{\mathbf{U}}^R(\xi = S_R, t) - \tilde{\mathbf{U}}^*(\xi = S_R, t)\right) =$$
$$\tilde{\mathbf{B}}(\mathbf{U}^*, \mathbf{U}^R)(\mathbf{U}^R - \mathbf{U}^*) + t\tilde{\mathbf{B}}(\mathbf{U}^*, \mathbf{U}^R)\left[\left(S_R\mathbf{I} - \mathbf{A}(\mathbf{U}^R)\right)(\partial_x \mathbf{U}^R) - \left(S_R\mathbf{I} - \mathbf{A}(\mathbf{U}^*)\right)(\partial_x \mathbf{U}^*)\right]$$
$$+ t\left\{\left[\frac{\partial \tilde{\mathbf{B}}(\mathbf{U}^*, \mathbf{U}^R)}{\partial (\mathbf{U}^*)_i}\right](\mathbf{U}^R - \mathbf{U}^*)\right\}\left[\left(S_R\mathbf{I} - \mathbf{A}(\mathbf{U}^*)\right)(\partial_x \mathbf{U}^*)\right]_i \quad (3.21)$$
$$+ t\left\{\left[\frac{\partial \tilde{\mathbf{B}}(\mathbf{U}^*, \mathbf{U}^R)}{\partial (\mathbf{U}^R)_i}\right](\mathbf{U}^R - \mathbf{U}^*)\right\}\left[\left(S_R\mathbf{I} - \mathbf{A}(\mathbf{U}^*)\right)(\partial_x \mathbf{U}^*)\right]_i$$

The repeated index "*i*" in the above equation denotes an Einstein summation over the components of the relevant vectors. The tensorial expressions in the last two terms of eqn. (3.21) are now clearly visible. It may be possible to evaluate them using implicit differentiation (i.e. Fréchêt derivatives). It is these tensorial terms that make the exact treatment of the GRP difficult for hyperbolic systems with non-conservative products. Using eqn. (3.21) and grouping all the *t*-dependent terms in eqn. (3.20) we get

$$\left(S_R\mathbf{I} - \mathbf{C}(\mathbf{U}^*)\right)\left(S_R\mathbf{I} - \mathbf{A}(\mathbf{U}^*)\right)(\partial_x \mathbf{U}^*) + S_R'(\mathbf{U}^* - \mathbf{U}^R) - \tilde{\mathbf{B}}(\mathbf{U}^*, \mathbf{U}^R)\left(S_R\mathbf{I} - \mathbf{A}(\mathbf{U}^*)\right)(\partial_x \mathbf{U}^*)$$
$$+ \left\{\left[\frac{\partial \tilde{\mathbf{B}}(\mathbf{U}^*, \mathbf{U}^R)}{\partial (\mathbf{U}^*)_i}\right](\mathbf{U}^R - \mathbf{U}^*)\right\}\left[\left(S_R\mathbf{I} - \mathbf{A}(\mathbf{U}^*)\right)(\partial_x \mathbf{U}^*)\right]_i =$$
$$\left(S_R\mathbf{I} - \mathbf{C}(\mathbf{U}^R)\right)\left(S_R\mathbf{I} - \mathbf{A}(\mathbf{U}^R)\right)(\partial_x \mathbf{U}^R) - \tilde{\mathbf{B}}(\mathbf{U}^*, \mathbf{U}^R)\left(S_R\mathbf{I} - \mathbf{A}(\mathbf{U}^R)\right)(\partial_x \mathbf{U}^R) \quad (3.22)$$
$$- \left\{\left[\frac{\partial \tilde{\mathbf{B}}(\mathbf{U}^*, \mathbf{U}^R)}{\partial (\mathbf{U}^R)_i}\right](\mathbf{U}^R - \mathbf{U}^*)\right\}\left[\left(S_R\mathbf{I} - \mathbf{A}(\mathbf{U}^*)\right)(\partial_x \mathbf{U}^*)\right]_i$$

Eqn. (3.22) is the general jump condition at second order that relates gradients in the right state $(\partial_x \mathbf{U}^R)$ to gradients in the resolved state $(\partial_x \mathbf{U}^*)$. We see that the presence of the tensorial terms make eqn. (3.22) very complicated. Except for the tensorial terms, eqn. (3.22) is not very complicated and it is indeed comparable in its complexity to eqn. (2.16).

It is possible to justify dropping the tensorial terms in eqn. (3.22) because they represent higher order derivatives in the non-conservative products. The current trend in the research literature is to write as much of the hyperbolic equation in conservation form as possible and to then deal with only a few terms that might be non-conservative. This provides an additional justification for dropping the tensorial terms in eqn. (3.22). The remaining part of eqn. (3.22) can be written compactly and elegantly as follows

$$\left[\left(S_R \mathbf{I} - \mathbf{C}(\mathbf{U}^*) - \tilde{\mathbf{B}}(\mathbf{U}^*, \mathbf{U}^R)\right)\left(S_R \mathbf{I} - \mathbf{A}(\mathbf{U}^*)\right)\right](\partial_x \mathbf{U}^*) + S_R'(\mathbf{U}^* - \mathbf{U}^R) =$$
$$\left[\left(S_R \mathbf{I} - \mathbf{C}(\mathbf{U}^R) - \tilde{\mathbf{B}}(\mathbf{U}^*, \mathbf{U}^R)\right)\left(S_R \mathbf{I} - \mathbf{A}(\mathbf{U}^R)\right)\right](\partial_x \mathbf{U}^R) \quad (3.23)$$

Comparing eqn. (3.23) to eqn. (2.16) we see an exact concordance. If the non-conservative terms could have been written as a Jacobian of some form of flux terms then eqn. (3.23) would have reduced exactly to eqn. (2.16). Of course, in reality this can't be done, but it helps to establish the concordance. Just as eqn. (3.16) has given us eqn. (3.23) at the right-going shock, we now use eqn. (3.17) at the left-going shock to get

$$\left[\left(S_L \mathbf{I} - \mathbf{C}(\mathbf{U}^*) - \tilde{\mathbf{B}}(\mathbf{U}^L, \mathbf{U}^*)\right)\left(S_L \mathbf{I} - \mathbf{A}(\mathbf{U}^*)\right)\right](\partial_x \mathbf{U}^*) + S_L'(\mathbf{U}^* - \mathbf{U}^L) =$$
$$\left[\left(S_L \mathbf{I} - \mathbf{C}(\mathbf{U}^L) - \tilde{\mathbf{B}}(\mathbf{U}^L, \mathbf{U}^*)\right)\left(S_L \mathbf{I} - \mathbf{A}(\mathbf{U}^L)\right)\right](\partial_x \mathbf{U}^L) \quad (3.24)$$

Eqns. (3.23) and (3.24) are extremely important for the GRP because they give us $(\partial_x \mathbf{U}^*)$ which leads us to the solution of the GRP for non-conservative hyperbolic systems. Eqns. (3.23) and (3.24) for hyperbolic systems with non-conservative products are analogous to eqns. (2.16) and (2.17) for conservation laws. The unknowns in the above two equations consist of the $M$-component vector $(\partial_x \mathbf{U}^*)$ and the two scalars $S_R'$ and $S_L'$. Consequently, we have $(M+2)$ unknowns whereas eqns. (3.23) and (3.24) constitute $2M$ linear equations in those unknowns. Therefore, for $M \geq 2$, a solution can always be found by minimizing the variation in eqns. (3.23) and (3.24) in a least squares sense. This completes the description of our solution methodology for obtaining the $M$-component gradient vector in the resolved state, $(\partial_x \mathbf{U}^*)$, as well as the curvature terms $S_R'$ and $S_L'$. By obtaining the above-mentioned $(M+2)$ unknowns we solve the HLL-GRP for hyperbolic systems with non-conservative products.

If one does not care to retain the curvature terms in the extremal waves, eqns. (3.23) and (3.24) can also be viewed as $2M$ equations in $M$ unknowns. This was the choice made in GBD. We have experimented with that possibility and found it to work very well even for hyperbolic systems with non-conservative products. For that reason, we used this approximation in all our numerical experiments.

## IV) HLLI-GRP-Based Scheme for Hyperbolic Systems

The easiest way to derive the fluctuation form of the HLL-GRP is to start with the conservation law given by eqn. (2.1) and cast it into fluctuation form. For this reason, Sub-section IV.a describes the use of the HLLI-GRP in the solution of a one-dimensional hyperbolic conservation law. Sub-section IV.b describes the formulation and use of the HLLI-GRP in the solution of a hyperbolic system in non-conservative form.

## IV.a) HLLI-GRP-Based Scheme for Hyperbolic Conservation Law

Consider the second order accurate solution of the hyperbolic conservation law given by eqn. (2.1) on a uniform mesh with zone size $\Delta x$ and a timestep $\Delta t$ using a GRP approach. Consider the update of a zone labeled "$j$" from a time $t^n$ to a time $t^{n+1} = t^n + \Delta t$. Within the zone "$j$" at time $t^n$ we have the mean state $\mathbf{U}_j^n$ and its gradient $(\partial_x \mathbf{U}_j^n)$. The gradient can be obtained by any monotonicity preserving limiter. We wish to evolve the solution vector to a time $t^{n+1}$ where the solution vector at the later time is denoted by $\mathbf{U}_j^{n+1}$. Say that the HLL Riemann problem at zone boundary $j-1/2$ at time $t^n$ has given us a spatially second order, but temporally first order, resolved state $\mathbf{U}_{HLL;j-1/2}^{*n}$ and entropy-satisfying flux $\mathbf{F}_{HLL;j-1/2}^{*n}$, see eqns. (2.8) and (2.9). A least squares solution of the GRP eqns. (2.16) and (2.17) at zone boundary $j-1/2$ at time $t^n$ gives us the gradient of the resolved state, $(\partial_x \mathbf{U}_{HLL;j-1/2}^{*n})$. At the zone boundary $j+1/2$ at time $t^n$ we obtain the analogous resolved state $\mathbf{U}_{HLL;j+1/2}^{*n}$, entropy-satisfying flux $\mathbf{F}_{HLL;j+1/2}^{*n}$, and gradient of the resolved state, $(\partial_x \mathbf{U}_{HLL;j+1/2}^{*n})$. Once we have the gradient of the resolved state, eqn. (2.21) shows us how to obtain the second order in space and time fluxes at each zone boundary. We can, therefore, write the one-step time update based on the GRP solver as

$$\mathbf{U}_j^{n+1} = \mathbf{U}_j^n - \frac{\Delta t}{\Delta x} \left\{ \begin{bmatrix} \mathbf{F}_{HLL;j+1/2}^{*n} - \frac{\Delta t}{2} \left( \mathbf{A}\left(\mathbf{U}_{HLL;j+1/2}^{*n}\right) \right)^2 \left( \partial_x \mathbf{U}_{HLL;j+1/2}^{*n} \right) \end{bmatrix} \\ - \begin{bmatrix} \mathbf{F}_{HLL;j-1/2}^{*n} - \frac{\Delta t}{2} \left( \mathbf{A}\left(\mathbf{U}_{HLL;j-1/2}^{*n}\right) \right)^2 \left( \partial_x \mathbf{U}_{HLL;j-1/2}^{*n} \right) \end{bmatrix} \right\} \quad (4.1)$$

If we are only interested in the update of a conservation law, the above equation is all we need. We can also rewrite the above equation as

$$\mathbf{U}_j^{n+1} = \mathbf{U}_j^n - \frac{\Delta t}{\Delta x} \left\{ \mathbf{F}_{HLL-GRP;j+1/2}^{n+1/2} - \mathbf{F}_{HLL-GRP;j-1/2}^{n+1/2} \right\} \quad (4.2)$$

With the time-centered, spatially and temporally second order fluxes given by

$$\mathbf{F}_{HLL-GRP;j+1/2}^{n+1/2} = \mathbf{F}_{HLL;j+1/2}^{*n} - \frac{\Delta t}{2} \left( \mathbf{A}\left(\mathbf{U}_{HLL;j+1/2}^{*n}\right) \right)^2 \left( \partial_x \mathbf{U}_{HLL;j+1/2}^{*n} \right) \quad (4.3)$$

We make the transcription $j+1/2 \to j-1/2$ in the above formula to get $\mathbf{F}_{HLL-GRP;j-1/2}^{n+1/2}$.

The previous expressions give us the fluxes explicitly for the HLL-GRP Riemann solver. Extending these fluxes to include the HLLI-GRP Riemann solver simply requires us to add the anti-diffusive flux contributions from eqns. (2.22) and (2.23). We get

$$\mathbf{F}_{HLLI-GRP;j+1/2}^{n+1/2} = \mathbf{F}_{HLL-GRP;j+1/2}^{*n+1/2} - \phi_{j+1/2} \frac{S_{R;j+1/2} S_{L;j+1/2}}{\left(S_{R;j+1/2} - S_{L;j+1/2}\right)} \sum_{p=1}^{P} \delta_{j+1/2}^{p} \left[ l_{j+1/2}^{p} \cdot \left(\mathbf{U}_{j+1/2}^{R;n+1/2} - \mathbf{U}_{j+1/2}^{L;n+1/2}\right) \right] r_{j+1/2}^{p}$$

(4.4)

We make the transcription $j+1/2 \rightarrow j-1/2$ in the above formula to get $\mathbf{F}_{HLLI-GRP;j-1/2}^{n+1/2}$. Please note the use of eqns. (2.18) and (2.19) in evaluating time-centered terms for the anti-diffusive flux contributions. We therefore have

$$\mathbf{U}_{j+1/2}^{R;n+1/2} = \mathbf{U}_{j+1}^{n} - \frac{1}{2}\left(\partial_x \mathbf{U}_{j+1}^{n}\right)\Delta x - \frac{\Delta t}{2}\mathbf{A}\left(\mathbf{U}_{j+1}^{n}\right)\left(\partial_x \mathbf{U}_{j+1}^{n}\right)$$

(4.5)

and

$$\mathbf{U}_{j+1/2}^{L;n+1/2} = \mathbf{U}_{j}^{n} + \frac{1}{2}\left(\partial_x \mathbf{U}_{j}^{n}\right)\Delta x - \frac{\Delta t}{2}\mathbf{A}\left(\mathbf{U}_{j}^{n}\right)\left(\partial_x \mathbf{U}_{j}^{n}\right)$$

(4.6)

The left eigenvectors $l_{j+1/2}^{p}$, and the right eigenvectors $r_{j+1/2}^{p}$, that are needed in eqn. (4.4), can also be evaluated by using eqn. (2.20) to obtain a time-centered resolved state in the Riemann fan given by

$$\mathbf{U}_{j+1/2}^{*;n+1/2} = \mathbf{U}_{j+1/2}^{*n} - \frac{\Delta t}{2}\mathbf{A}\left(\mathbf{U}_{j+1/2}^{*n}\right)\left(\partial_x \mathbf{U}_{j+1/2}^{*n}\right)$$

(4.7)

The use of $\phi_{j+1/2}$ as a flattener function is described in Colella and Woodward (1984) or Balsara (2012). The expression for $\delta_{j+1/2}^{p}$ is still given by eqn. (2.23). Once the fluxes from the HLLI-GRP are in hand, the update equation for the scheme, i.e. eqn. (4.2), can be re-written with the transcription $HLL-GRP \rightarrow HLLI-GRP$. This completes our description of how the HLL-GRP and HLLI-GRP are to be used to obtain schemes for conservation laws.

### IV.b) HLLI-GRP-Based Scheme for Non-Conservative Hyperbolic Systems

For a hyperbolic system in non-conservative form, we necessarily need to use a fluctuation form for the time update of the PDE. However, we can start from the conservation form and recast it in fluctuation form. This gives us a form for the time update that is equally suitable for hyperbolic systems that are in conservation and non-conservative form.

Realize therefore that we are also interested in deriving an update strategy by modifying eqn. (4.1) so that it is written in fluctuation form. Our first step is the suitable insertion of the flux term $\mathbf{F}(\mathbf{U}_j^n)$ in eqn. (4.1) so as to recast that equation in a form that is closer to the fluctuation form. We therefore rewrite eqn. (4.1) in a form that more closely mimics the fluctuation form as follows

$$\mathbf{U}_j^{n+1} = \mathbf{U}_j^n - \frac{\Delta t}{\Delta x}\left\{\begin{aligned}&\left[\mathbf{F}_{HLL;j+1/2}^{*n} - \mathbf{F}(\mathbf{U}_j^n) - \frac{\Delta t}{2}\left(\mathbf{A}(\mathbf{U}_{HLL;j+1/2}^{*n})\right)^2\left(\partial_x \mathbf{U}_{HLL;j+1/2}^{*n}\right)\right]\\&+\left[\mathbf{F}(\mathbf{U}_j^n) - \mathbf{F}_{HLL;j-1/2}^{*n} + \frac{\Delta t}{2}\left(\mathbf{A}(\mathbf{U}_{HLL;j-1/2}^{*n})\right)^2\left(\partial_x \mathbf{U}_{HLL;j-1/2}^{*n}\right)\right]\end{aligned}\right\} \quad (4.8)$$

The above equation is still an intermediate equation. It has some ingredients of a fluctuation form, but it is not truly in fluctuation form. This is because it is still not possible to identify $\mathbf{F}_{HLL;j+1/2}^{*n} - \mathbf{F}(\mathbf{U}_j^n)$ with the left-going fluctuation $\mathbf{D}_{HLL;j+1/2}^{-;n}$ that is evaluated at zone boundary $j+1/2$ at a time $t^n$. Similarly, it is still not possible to identify $\mathbf{F}(\mathbf{U}_j^n) - \mathbf{F}_{HLL;j-1/2}^{*n}$ with the right-going fluctuation $\mathbf{D}_{HLL;j-1/2}^{+;n}$ that is evaluated at zone boundary $j-1/2$ at a time $t^n$. The correct identification of the left-going fluctuation $\mathbf{D}_{HLL;j+1/2}^{-;n}$ that is evaluated at zone boundary $j+1/2$ at a time $t^n$ is given by

$$\mathbf{D}_{HLL;j+1/2}^{-;n} = \mathbf{F}_{HLL;j+1/2}^{*n} - \mathbf{F}\left(\mathbf{U}_j^n + \frac{1}{2}(\partial_x \mathbf{U}_j^n)\Delta x\right) \cong \mathbf{F}_{HLL;j+1/2}^{*n} - \mathbf{F}(\mathbf{U}_j^n) - \frac{1}{2}\mathbf{A}(\mathbf{U}_j^n)(\partial_x \mathbf{U}_j^n)\Delta x \quad (4.9a)$$

Likewise, the correct identification of the right-going fluctuation $\mathbf{D}_{HLL;j-1/2}^{+;n}$ that is evaluated at zone boundary $j-1/2$ at a time $t^n$ is given by

$$\mathbf{D}_{HLL;j-1/2}^{+;n} = \mathbf{F}\left(\mathbf{U}_j^n - \frac{1}{2}(\partial_x \mathbf{U}_j^n)\Delta x\right) - \mathbf{F}_{HLL;j-1/2}^{*n} \cong \mathbf{F}(\mathbf{U}_j^n) - \mathbf{F}_{HLL;j-1/2}^{*n} - \frac{1}{2}\mathbf{A}(\mathbf{U}_j^n)(\partial_x \mathbf{U}_j^n)\Delta x \quad (4.9b)$$

The above two equations can be written in a form that is more useful for incorporation in eqn. (4.8) as follows

$$\mathbf{F}_{HLL;j+1/2}^{*n} - \mathbf{F}(\mathbf{U}_j^n) = \mathbf{D}_{HLL;j+1/2}^{-;n} + \frac{1}{2}\mathbf{A}(\mathbf{U}_j^n)(\partial_x \mathbf{U}_j^n)\Delta x \quad (4.11)$$

and

$$\mathbf{F}(\mathbf{U}_j^n) - \mathbf{F}_{HLL;j-1/2}^{*n} = \mathbf{D}_{HLL;j-1/2}^{+;n} + \frac{1}{2}\mathbf{A}(\mathbf{U}_j^n)(\partial_x \mathbf{U}_j^n)\Delta x \quad (4.12)$$

Replacing the above two equations into eqn. (4.8) gives

$$\mathbf{U}_j^{n+1} = \mathbf{U}_j^n - \Delta t \mathbf{A}\left(\mathbf{U}_j^n\right)\left(\partial_x \mathbf{U}_j^n\right) - \frac{\Delta t}{\Delta x} \left\{ \begin{aligned} &\left[ \mathbf{D}_{HLL;j+1/2}^{-;n} - \frac{\Delta t}{2}\left(\mathbf{A}\left(\mathbf{U}_{HLL;j+1/2}^{*n}\right)\right)^2 \left(\partial_x \mathbf{U}_{HLL;j+1/2}^{*n}\right) \right] \\ &+ \left[ \mathbf{D}_{HLL;j-1/2}^{+;n} + \frac{\Delta t}{2}\left(\mathbf{A}\left(\mathbf{U}_{HLL;j-1/2}^{*n}\right)\right)^2 \left(\partial_x \mathbf{U}_{HLL;j-1/2}^{*n}\right) \right] \end{aligned} \right\} \quad (4.13)$$

The above equation can be most compactly written as

$$\mathbf{U}_j^{n+1} = \mathbf{U}_j^n - \Delta t \mathbf{A}\left(\mathbf{U}_j^n\right)\left(\partial_x \mathbf{U}_j^n\right) - \frac{\Delta t}{\Delta x}\left\{ \mathbf{D}_{HLL-GRP;j+1/2}^{-;n+1/2} + \mathbf{D}_{HLL-GRP;j-1/2}^{+;n+1/2} \right\} \quad (4.14)$$

With the additional definitions

$$\mathbf{D}_{HLL-GRP;j+1/2}^{-;n+1/2} = \mathbf{D}_{HLL;j+1/2}^{-;n} - \frac{\Delta t}{2}\left(\mathbf{A}\left(\mathbf{U}_{HLL;j+1/2}^{*n}\right)\right)^2 \left(\partial_x \mathbf{U}_{HLL;j+1/2}^{*n}\right) \quad (4.15)$$

and

$$\mathbf{D}_{HLL-GRP;j-1/2}^{+;n+1/2} = \mathbf{D}_{HLL;j-1/2}^{+;n} + \frac{\Delta t}{2}\left(\mathbf{A}\left(\mathbf{U}_{HLL;j-1/2}^{*n}\right)\right)^2 \left(\partial_x \mathbf{U}_{HLL;j-1/2}^{*n}\right) \quad (4.16)$$

It is now evident, that the above two equations give us the time-centered, spatially and temporally second order fluctuations for the HLL-GRP. Also please notice that the left-going fluctuation $\mathbf{D}_{HLL;j+1/2}^{-;n}$ and the right going fluctuation $\mathbf{D}_{HLL;j+1/2}^{+;n}$ at the zone boundary $j+1/2$ will satisfy the consistency condition for the fluctuations given by eqn. (3.12) at time $t^n$. This consistency is guaranteed by the construction of the HLL Riemann solver for hyperbolic systems with non-conservative products. The structure of the above two equations is such that the time-dependence also guarantees that the time-centered fluctuations $\mathbf{D}_{HLL-GRP;j+1/2}^{+;n+1/2}$ and $\mathbf{D}_{HLL-GRP;j+1/2}^{-;n+1/2}$ will also satisfy the same consistency condition for the fluctuations at the centered time $t^{n+1/2}$. Please also realize that it is not profitable to use a fluctuation form when solving for a hyperbolic system in conservation form. However, eqns. (4.14), (4.15) and (4.16) are very useful when solving a hyperbolic system with non-conservative products.

When dealing with a hyperbolic system that has non-conservative products, we iterate eqn. (3.9) to convergence in order to obtain the resolved state in the HLL Riemann solver. The spatially second order accurate fluctuations at time $t^n$ at each zone boundary can then be obtained from eqns. (3.10) and (3.11). Specifically, we can now write

$$\mathbf{D}_{HLL;j+1/2}^{-;n} = S_{L;j+1/2}\left( \mathbf{U}_{HLL;j+1/2}^{*n} - \mathbf{U}_j^n - \frac{1}{2}\left(\partial_x \mathbf{U}_j^n\right)\Delta x \right) \quad (4.17)$$

and

$$\mathbf{D}^{+;n}_{HLL;j-1/2} = S_{R;j-1/2}\left(\mathbf{U}^n_j - \frac{1}{2}(\partial_x\mathbf{U}^n_j)\Delta x - \mathbf{U}^{*n}_{HLL;j-1/2}\right) \tag{4.18}$$

Eqns. (3.23) and (3.24) (which are analogous to eqns. (2.16) and (2.17)) can then be used to obtain the spatial gradient of the resolved state in the HLL Riemann solver. Once the resolved state and its spatial gradient are in hand at each zone boundary, eqns. (4.15) and (4.16) give us the spatially and temporally second order fluctuations at each zone boundary for the HLL-GRP as it applies to non-conservative hyperbolic systems.

Obtaining the spatially and temporally second order fluctuations for the HLLI-GRP is now easy. We just have to provide the anti-diffusive contributions from eqns. (3.13) and (3.14) to the HLL-GRP from eqns. (4.15) and (4.16). However, notice that the anti-diffusive contributions from eqns. (3.13) and (3.14) have to be slightly altered because we want a time-centered contribution. This is easily accomplished by making the formal transcription $(\mathbf{U}^R - \mathbf{U}^L) \to (\mathbf{U}^{R;n+1/2} - \mathbf{U}^{L;n+1/2})$ in eqns. (3.13) and (3.14) for the anti-diffusive contributions to the fluctuations. Here $\mathbf{U}^{R;n+1/2}_{j+1/2}$ and $\mathbf{U}^{L;n+1/2}_{j+1/2}$ are obtained from eqns. (4.5) and (4.6). We therefore assemble our final expressions for the time-centered fluctuations for the HLLI-GRP as

$$\mathbf{D}^{-;n+1/2}_{HLLI-GRP;j+1/2} = \mathbf{D}^{-;n+1/2}_{HLL-GRP;j+1/2} - \phi_{j+1/2}\frac{S_{R;j+1/2}S_{L;j+1/2}}{(S_{R;j+1/2} - S_{L;j+1/2})}\sum_{p=1}^{P}\delta^p_{j+1/2}\left[l^p_{j+1/2}\cdot\left(\mathbf{U}^{R;n+1/2}_{j+1/2} - \mathbf{U}^{L;n+1/2}_{j+1/2}\right)\right]r^p_{j+1/2}$$
(4.19)

and

$$\mathbf{D}^{+;n+1/2}_{HLLI-GRP;j-1/2} = \mathbf{D}^{+;n+1/2}_{HLL-GRP;j-1/2} + \phi_{j-1/2}\frac{S_{R;j-1/2}S_{L;j-1/2}}{(S_{R;j-1/2} - S_{L;j-1/2})}\sum_{p=1}^{P}\delta^p_{j-1/2}\left[l^p_{j-1/2}\cdot\left(\mathbf{U}^{R;n+1/2}_{j-1/2} - \mathbf{U}^{L;n+1/2}_{j-1/2}\right)\right]r^p_{j-1/2}$$
(4.20)

As before, the left eigenvectors $l^p_{j+1/2}$, and the right eigenvectors $r^p_{j+1/2}$, are evaluated using the time-centered resolved state in eqn. (4.7). Please also notice that the structure of the above two equations ensures that the time-centered fluctuations $\mathbf{D}^{+;n+1/2}_{HLLI-GRP;j+1/2}$ and $\mathbf{D}^{-;n+1/2}_{HLLI-GRP;j+1/2}$ will also satisfy the same consistency condition for the fluctuations at the centered time $t^{n+1/2}$. Once the fluctuations from the HLLI-GRP are in hand, the update equation for the scheme, i.e. eqn. (4.14), can be re-written with the transcription $HLL-GRP \to HLLI-GRP$. This completes our description of the scheme that uses the HLLI-GRP for hyperbolic systems with non-conservative products.

**V) Inclusion of Stiff Source Terms in the GRP Solver**

Several physical problems involving hyperbolic systems also require the inclusion of stiff source terms. In Sub-section V.a we provide a quick description of a strategy for including stiff source terms. In Sub-section V.b we show how this strategy is incorporated for conservation laws involving stiff source terms. In Sub-section V.c we show how this strategy is incorporated for hyperbolic systems with non-conservative products and stiff source terms.

**V.a) Quick Description of ADER Scheme for Stiff Source Terms**

Thus far, we have been concerned with the GRP solver for the $M \times M$ conservation law given by eqn. (2.1) or the analogous hyperbolic system with non-conservative products given by eqn. (3.1). There are no source terms in those equations. For certain hyperbolic problems we may have source terms so that eqn. (3.1) becomes

$$\frac{\partial \mathbf{U}}{\partial t} + \frac{\partial \mathbf{F}(\mathbf{U})}{\partial x} + \mathbf{B}(\mathbf{U})\frac{\partial \mathbf{U}}{\partial x} = \mathbf{S}(\mathbf{U}) \tag{5.1}$$

Just as in eqn. (3.1), "$\mathbf{U}$" is an $M$-component vector, "$\mathbf{F}(\mathbf{U})$" is an $M$-component flux that depends on "$\mathbf{U}$" and the $M \times M$ matrix "$\mathbf{B}(\mathbf{U})$" contains the non-conservative products. The new feature is that $\mathbf{S}(\mathbf{U})$ is an $M$-component vector of potentially stiff source terms that depend on "$\mathbf{U}$". Since the case with stiff source terms is more interesting, and more challenging, we direct out attention to that case. The definition of the $M \times M$ characteristic matrix "$\mathbf{A}(\mathbf{U})$" is unchanged so that we can write eqn. (5.1) as

$$\frac{\partial \mathbf{U}}{\partial t} + \mathbf{A}(\mathbf{U})\frac{\partial \mathbf{U}}{\partial x} = \mathbf{S}(\mathbf{U}) \tag{5.2}$$

Whether inside a zone, or inside a Riemann fan, we desire a solution strategy that starts with initial conditions given by a solution vector "$\mathbf{U}^0$" and its spatial gradient "$(\partial_x \mathbf{U}^0)$" at time $t^n$ and produces an evolution-in-the-small in space and time that is consistent with the governing equations eqn. (5.2). In other words, we seek a solution in space and time that is valid up to second order in the spatial distance $\Delta x$ and the temporal distance $\Delta t$. (Please also note that the initial gradient is re-defined in this Section to be $(\partial_x \mathbf{U}^0)$, instead of $(\partial_x \mathbf{U}^1)$ as in eqn. (2.3), because this choice makes the notation in this section more intuitive.) The characteristic matrix can be treated explicitly in the solution that we seek, however, the source terms have to be evolved implicitly when they are stiff. The ADER scheme is designed to give us exactly such a solution strategy. In Balsara *et al.* ([7]; henceforth B17) we designed an ultra-efficient second-order accurate ADER scheme for stiff source terms that accomplishes this.

We transcribe that ADER scheme from B17 to meet our needs in this Sub-section. When the source terms are stiff, they have a strong influence not just on the temporal evolution of the solution vector but also on the temporal evolution of its gradient. Consequently, at a time $t^n + \Delta t/2$ we seek a solution vector "$\mathbf{U}^{1/2}$" and its spatial gradient "$(\partial_x \mathbf{U}^{1/2})$", furthermore, at a time of $t^n + \Delta t$ seek a solution vector "$\mathbf{U}^1$" and its spatial gradient "$(\partial_x \mathbf{U}^1)$". Transcribing eqn. (4.3) of B17, and with the initial solution centered at the origin $(x,t) = (0,0)$, we seek a space-time expansion of the form

$$\mathbf{U}(x,t) = \left\{\mathbf{U}^{1/2} + (\partial_x \mathbf{U}^{1/2})x\right\}\left(2 - 2\frac{t}{\Delta t}\right) + \left\{\mathbf{U}^1 + (\partial_x \mathbf{U}^1)x\right\}\left(2\frac{t}{\Delta t} - 1\right) \tag{5.3}$$

This expansion is linear in space and time, thereby being second order accurate. Moreover, by retaining the time evolution of the gradients, it permits stiff source terms to influence not just the solution but also the gradient of the solution. If the source terms are very stiff, they can even make large changes to the solution in a fraction of a timestep; and by positing the modes of our solution at times $t^n + \Delta t/2$ and $t^n + \Delta t$, we allow for that to happen. These are all the desirable traits that we want in our solution. The above expansion has to be solved for in a fashion that is consistent with the initial conditions and the governing equation. The initial conditions can themselves be expressed as

$$\mathbf{W}(x) = \mathbf{U}^0 + (\partial_x \mathbf{U}^0)x \tag{5.4}$$

Moreover, it should treat the source terms in a time-implicit fashion. The source terms are given by an expansion that is similar to eqn. (5.3) and that expansion is given by

$$\begin{aligned}\mathbf{S}(x,t) &= \left\{\mathbf{S}(\mathbf{U}^{1/2}) + [\partial_\mathbf{U}\mathbf{S}(\mathbf{U}^{1/2})](\partial_x \mathbf{U}^{1/2})x\right\}\left(2 - 2\frac{t}{\Delta t}\right) \\ &+ \left\{\mathbf{S}(\mathbf{U}^1) + [\partial_\mathbf{U}\mathbf{S}(\mathbf{U}^1)](\partial_x \mathbf{U}^1)x\right\}\left(2\frac{t}{\Delta t} - 1\right)\end{aligned} \tag{5.5}$$

Notice that eqns. (5.3) and (5.5) have a similar structure. This is because, as the solution undergoes large changes in response to stiff source terms, we also want to make allowance for the source terms to undergo substantial, self-consistent changes in response to the rapidly-changing solution.

Eqns. (5.3) and (5.5) are to be viewed in a Galerkin sense as a basis expansion in space and time. It is easy to identify the space-time trial functions in eqn. (5.3). By making the coefficients of the basis expansion in eqn. (5.3) different from the coefficients of the initial conditions in eqn. (5.4), we allow for the fact that the source terms could have a very strong influence on the initial conditions. In other words, we allow the source terms to be stiff and we allow for the fact that the source terms can make large changes to the solution within a time interval that can be less than the timestep $\Delta t$. Using test functions that are the same as the trial functions, we make a Galerkin projection of the governing eqn. (5.2). The projection is made so that it is consistent with the initial conditions in eqn. (5.4). Transcribing eqns. (4.10) and (4.11) of B17 we get

$$\mathbf{U}^{1/2} = \mathbf{U}^0 + \frac{4}{6}\Delta t\ \mathbf{S}(\mathbf{U}^{1/2}) - \frac{1}{6}\Delta t\ \mathbf{S}(\mathbf{U}^1) - \frac{4}{6}\Delta t\ \mathbf{A}(\mathbf{U}^{1/2})(\partial_x \mathbf{U}^{1/2}) + \frac{1}{6}\Delta t\ \mathbf{A}(\mathbf{U}^1)(\partial_x \mathbf{U}^1) \quad (5.6)$$

$$(\partial_x \mathbf{U}^{1/2}) = (\partial_x \mathbf{U}^0) + \frac{4}{6}\Delta t\ \left[\partial_\mathbf{U} \mathbf{S}(\mathbf{U}^{1/2})\right](\partial_x \mathbf{U}^{1/2}) - \frac{1}{6}\Delta t\ \left[\partial_\mathbf{U} \mathbf{S}(\mathbf{U}^1)\right](\partial_x \mathbf{U}^1) \quad (5.7)$$

Likewise, transcribing eqns. (4.14) and (4.15) of B17 we get

$$\mathbf{U}^1 = \mathbf{U}^0 + \Delta t\ \mathbf{S}(\mathbf{U}^{1/2}) - \Delta t\ \mathbf{A}(\mathbf{U}^{1/2})(\partial_x \mathbf{U}^{1/2}) \quad (5.8)$$

$$(\partial_x \mathbf{U}^1) = (\partial_x \mathbf{U}^0) + \Delta t\ \left[\partial_\mathbf{U} \mathbf{S}(\mathbf{U}^{1/2})\right](\partial_x \mathbf{U}^{1/2}) \quad (5.9)$$

Here "$\left[\partial_\mathbf{U} \mathbf{S}(\mathbf{U}^{1/2})\right]$" denotes the Jacobian of the source term with respect to the solution vector, evaluated at "$\mathbf{U} = \mathbf{U}^{1/2}$". Likewise, "$\left[\partial_\mathbf{U} \mathbf{S}(\mathbf{U}^1)\right]$" denotes the Jacobian of the source term with respect to the solution vector, evaluated at "$\mathbf{U} = \mathbf{U}^1$". This completes our description of the governing equations in our ADER scheme.

The system of equations given by eqns. (5.6) to (5.9) can be solved via fixed point (Picard) iteration with the initial guess given by $\mathbf{U}^1 = \mathbf{U}^{1/2} = \mathbf{U}^0$ and $(\partial_x \mathbf{U}^1) = (\partial_x \mathbf{U}^{1/2}) = (\partial_x \mathbf{U}^0)$. They converge to the desired accuracy within about two iterations. When the source terms are non-stiff, the source terms don't require implicit treatment. When the source terms are not extremely stiff, one can get by with treating just eqns. (5.6) and (5.8) implicitly while treating eqns. (5.7) and (5.9) as auxiliary equations. In that case, Section IV.3 of B17 shows that eqns. (5.6) and (5.8) can be solved via the solution of one $M \times M$ linear system and the inversion of one $M \times M$ matrix. In practice, we just introduce one implicit sub-iteration step for the source terms within the Picard iteration. This sub-iteration is an optimal solution strategy because we are getting a second-order, fully implicit treatment of source terms at the computational cost that is slightly less than the inversion of two $M \times M$ matrices. When the source terms are extremely stiff, the gradient terms in eqns. (5.7) and (5.9) should also be treated implicitly. In that case, Section IV.3 of B17 shows that we get an implicit system that is entirely analogous to the one formed by eqns. (5.6) and (5.8) and the same optimal solution strategy can be applied to the sub-iteration of the gradients. For further details, please see Section IV of B17.

This completes our description of the ADER scheme for treating stiff source terms, especially as it is modified for the GRP solver. A similar approach of using the ADER scheme to form a GRP solver that can handle stiff sources was tried in Montecinos and Toro [44] and Goetz and Dumbser [33]. In the next sub-section, we show how it is incorporated into the GRP solver in order to obtain a GRP solver-based numerical scheme that is stable in the presence of stiff source terms.

**V.b) HLLI-GRP-Based Scheme for Hyperbolic Conservation Law With Stiff Source Terms**

Consider the second order accurate solution of the hyperbolic conservation law given by eqn. (5.1) on a mesh with zone size $\Delta x$ and a timestep $\Delta t$ using a GRP approach. In this Sub-section we consider the situation where the non-conservative products are not present. All of the notation from Sub-sections IV.a and IV.b is taken over for this and the next Sub-section. Within the zone "$j$" at time $t^n$ we have the mean state $\mathbf{U}_j^n$ and its gradient $\left(\partial_x \mathbf{U}_j^n\right)$. The crucial difference is that the in-the-small time-evolution in that zone is now expressed by using eqn. (5.3) and can be written as

$$\mathbf{U}_j(x,t) = \left\{\mathbf{U}_j^{n+1/2} + \left(\partial_x \mathbf{U}_j^{n+1/2}\right)x\right\}\left(2 - 2\frac{t}{\Delta t}\right) + \left\{\mathbf{U}_j^{n+1} + \left(\partial_x \mathbf{U}_j^{n+1}\right)x\right\}\left(2\frac{t}{\Delta t} - 1\right) \tag{5.10}$$

The initial conditions can themselves be expressed as

$$\mathbf{W}(x) = \mathbf{U}_j^0 + \left(\partial_x \mathbf{U}_j^0\right)x \tag{5.11}$$

The corresponding expression for the source terms is given by

$$\begin{aligned}\mathbf{S}_j(x,t) &= \left\{\mathbf{S}\left(\mathbf{U}_j^{n+1/2}\right) + \left[\partial_{\mathbf{U}}\mathbf{S}\left(\mathbf{U}_j^{n+1/2}\right)\right]\left(\partial_x \mathbf{U}_j^{n+1/2}\right)x\right\}\left(2 - 2\frac{t}{\Delta t}\right) \\ &+ \left\{\mathbf{S}\left(\mathbf{U}_j^{n+1}\right) + \left[\partial_{\mathbf{U}}\mathbf{S}\left(\mathbf{U}_j^{n+1}\right)\right]\left(\partial_x \mathbf{U}_j^{n+1}\right)x\right\}\left(2\frac{t}{\Delta t} - 1\right)\end{aligned} \tag{5.12}$$

The terms $\mathbf{U}_j^{n+1/2}$, $\left(\partial_x \mathbf{U}_j^{n+1/2}\right)$, $\mathbf{U}_j^{n+1}$ and $\left(\partial_x \mathbf{U}_j^{n+1}\right)$ in eqn. (5.10) are obtained by iterating with the help of eqns. (5.6) to (5.9). This iteration is local and restricted to the zone "$j$"; and usually about two iterations are sufficient. The eqn. (5.12) for the source term will also have to be averaged over space and time to obtain its contribution to the update equation.

At the zone boundary $j+1/2$ we will still have the spatially second order, but temporally first order, resolved state $\mathbf{U}_{HLL;j+1/2}^{*n}$ from the Riemann solver. The source terms do not influence the least squares solution of the GRP eqns. (2.16) and (2.17) with the result that the gradient of the resolved state at the same zone boundary is still given by $\left(\partial_x \mathbf{U}_{HLL;j+1/2}^{*n}\right)$. Now notice something interesting about eqns. (5.6) to (5.9) – they only depend on the solution vector, its gradient and $\Delta t$; but they do not involve $\Delta x$. This has the interesting result that they are also applicable at the zone boundary. At the zone boundary $j+1/2$, they can be applied *within* the Riemann fan! As in the previous paragraph, we start with initial conditions given by $\mathbf{U}_{HLL;j+1/2}^{*n}$ and $\left(\partial_x \mathbf{U}_{HLL;j+1/2}^{*n}\right)$ within the Riemann fan. We use these initial conditions, along with eqns. (5.6) to (5.9), to obtain $\mathbf{U}_{HLL;j+1/2}^{*n+1/2}$

, $\left(\partial_x \mathbf{U}^{*n+1/2}_{HLL;j+1/2}\right)$, $\mathbf{U}^{*n+1}_{HLL;j+1/2}$ and $\left(\partial_x \mathbf{U}^{*n+1}_{HLL;j+1/2}\right)$. At the zone boundary, we are only interested in the temporal evolution so that we have the analogue of eqn. (5.10) given by

$$\mathbf{U}^*_{j+1/2}\left(x = x_{j+1/2}, t\right) = \mathbf{U}^{*n+1/2}_{HLL;j+1/2}\left(2 - 2\frac{t}{\Delta t}\right) + \mathbf{U}^{*n+1}_{HLL;j+1/2}\left(2\frac{t}{\Delta t} - 1\right) \tag{5.13}$$

We can now write the update that is analogous to eqn. (4.2) as follows

$$\mathbf{U}^{n+1}_j = \mathbf{U}^n_j - \frac{\Delta t}{\Delta x}\left\{\mathbf{F}^{n+1/2}_{HLL-GRP;j+1/2} - \mathbf{F}^{n+1/2}_{HLL-GRP;j-1/2}\right\} + \Delta t \; \mathbf{S}\left(\mathbf{U}^{n+1/2}_j\right) \tag{5.14}$$

With the time-centered, spatially and temporally second order fluxes given by

$$\mathbf{F}^{n+1/2}_{HLL-GRP;j+1/2} = \mathbf{F}^{*n}_{HLL;j+1/2} + \mathbf{A}\left(\mathbf{U}^{*n+1/2}_{HLL;j+1/2}\right)\left(\mathbf{U}^{*n+1}_{HLL;j+1/2} - \mathbf{U}^{*n+1/2}_{HLL;j+1/2}\right) \tag{5.15}$$

We have taken the liberty of making the transcription $\mathbf{A}\left(\mathbf{U}^{*n}_{HLL;j+1/2}\right) \to \mathbf{A}\left(\mathbf{U}^{*n+1/2}_{HLL;j+1/2}\right)$ in the above equation in order to obtain a more time-centered flux. We make the transcription $j+1/2 \to j-1/2$ in the above formula to get $\mathbf{F}^{n+1/2}_{HLL-GRP;j-1/2}$.

To make the extension to the HLLI-GRP solver, we keep eqn. (4.4) for the HLLI-GRP flux (with its anti-diffusive contribution) unchanged. However, we should use eqn. (5.10) to get $\mathbf{U}^{L;n+1/2}_{j+1/2}$ so that we have

$$\mathbf{U}^{L;n+1/2}_{j+1/2} = \mathbf{U}^{n+1/2}_j + \frac{\Delta x}{2}\left(\partial_x \mathbf{U}^{n+1/2}_j\right) \tag{5.16}$$

and

$$\mathbf{U}^{R;n+1/2}_{j+1/2} = \mathbf{U}^{n+1/2}_{j+1} - \frac{\Delta x}{2}\left(\partial_x \mathbf{U}^{n+1/2}_{j+1}\right) \tag{5.17}$$

The time-centered solution vector and its gradient already incorporate the influence of the stiff source term and the characteristic matrix on the time evolution. This is because the source terms and the characteristic matrix were involved in the iteration of eqns. (5.6) to (5.9). The left eigenvectors $l^p_{j+1/2}$, and the right eigenvectors $r^p_{j+1/2}$, that are needed in eqn. (4.4), can also be evaluated by using the state vector $\mathbf{U}^{*n+1/2}_{HLL;j+1/2}$ from eqn. (5.13). Again, the contributions from the stiff source terms and the characteristic matrix on the time evolution are already built into $\mathbf{U}^{*n+1/2}_{HLL;j+1/2}$. Once the fluxes from the HLLI-GRP are in hand, the update equation for the scheme, i.e. eqn. (5.14), can be re-written with the transcription $HLL-GRP \to HLLI-GRP$. This completes our description of the scheme that uses the HLLI-GRP solver for conservation laws with stiff source terms.

## V.c) HLLI-GRP-Based Scheme for Non-Conservative Hyperbolic Systems With Stiff Source Terms

When non-conservative products are present, the update equation that is analogous to eqn. (4.14) becomes

$$\mathbf{U}_j^{n+1} = \mathbf{U}_j^n - \Delta t \mathbf{A}\left(\mathbf{U}_j^n\right)\left(\partial_x \mathbf{U}_j^n\right) - \frac{\Delta t}{\Delta x}\left\{\mathbf{D}_{HLL-GRP;j+1/2}^{-;n+1/2} + \mathbf{D}_{HLL-GRP;j-1/2}^{+;n+1/2}\right\} + \Delta t\, \mathbf{S}\left(\mathbf{U}_j^{n+1/2}\right) \quad (5.18)$$

With the additional definitions

$$\mathbf{D}_{HLL-GRP;j+1/2}^{-;n+1/2} = \mathbf{D}_{HLL;j+1/2}^{-;n} + \mathbf{A}\left(\mathbf{U}_{HLL;j+1/2}^{*n+1/2}\right)\left(\mathbf{U}_{HLL;j+1/2}^{*n+1} - \mathbf{U}_{HLL;j+1/2}^{*n+1/2}\right) \quad (5.19)$$

and

$$\mathbf{D}_{HLL-GRP;j-1/2}^{+;n+1/2} = \mathbf{D}_{HLL;j-1/2}^{+;n} - \mathbf{A}\left(\mathbf{U}_{HLL;j-1/2}^{*n+1/2}\right)\left(\mathbf{U}_{HLL;j-1/2}^{*n+1} - \mathbf{U}_{HLL;j-1/2}^{*n+1/2}\right) \quad (5.20)$$

Eqns. (4.17) and (4.18) remain unchanged for $\mathbf{D}_{HLL;j+1/2}^{-;n}$ and $\mathbf{D}_{HLL;j-1/2}^{+;n}$ respectively. In the above two equations we have taken the liberty of making the transcription $\mathbf{A}\left(\mathbf{U}_{HLL;j+1/2}^{*n}\right) \rightarrow \mathbf{A}\left(\mathbf{U}_{HLL;j+1/2}^{*n+1/2}\right)$ in order to obtain more time-centered expressions for the fluctuations.

Eqns. (4.19) and (4.20) for the anti-diffusive contributions to the fluctuations remain unchanged as long as eqns. (5.16) and (5.17) are used for $\mathbf{U}_{j+1/2}^{L;n+1/2}$ and $\mathbf{U}_{j+1/2}^{R;n+1/2}$. The left eigenvectors $l_{j+1/2}^p$, and the right eigenvectors $r_{j+1/2}^p$, that are needed in eqns. (4.19) and (4.20), can also be evaluated by using the state vector $\mathbf{U}_{HLL;j+1/2}^{*n+1/2}$ from eqn. (5.13). Again, the contributions from the stiff source terms to the time evolution are already built into $\mathbf{U}_{HLL;j+1/2}^{*n+1/2}$. Once the fluctuations from the HLLI-GRP are in hand, the update equation for the scheme, i.e. eqn. (5.18), can be re-written with the transcription $HLL-GRP \rightarrow HLLI-GRP$. This completes our description of the scheme that uses the HLLI-GRP solver for hyperbolic systems with non-conservative products and stiff source terms.

## VI) Test Problems

To illustrate the versatility of our HLLI-GRP we show that it works well for several different hyperbolic systems of practical interest. The hyperbolic systems we pick include the Euler equations, the MHD equations, the shallow water equations with spatially varying bathymetry and the compressible Navier Stokes equations. The Euler and MHD equations provide examples of a hyperbolic system in conservation form. The shallow water equations serve as an example of a non-conservative system. When the bathymetry is non-constant, we have an example of a non-conservative system with stiff source terms. The compressible Navier Stokes equations in relaxation form give us an example of a conservation law with stiff source terms. As a result, our examples cover all the different types of situations considered in this paper.

All test problems were run with a second order code accurate that used piecewise linear reconstruction using an MC limiter. The HLLI-GRP was then used to provide the fluxes or fluctuations depending on the system being considered. All tests were run with a CFL of 0.8.

## VI.a) Euler Equations

The Euler equations in one dimension are given by

$$\frac{\partial}{\partial t}\begin{pmatrix} \rho \\ \rho v_x \\ \rho v_y \\ \rho v_z \\ \varepsilon \end{pmatrix} + \frac{\partial}{\partial x}\begin{pmatrix} \rho v_x \\ \rho v_x^2 + P \\ \rho v_x v_y \\ \rho v_x v_z \\ (\varepsilon + P) v_x \end{pmatrix} = 0 \qquad (6.1)$$

Here $\rho$ is the gas density, "P" is the gas pressure and $v_x$, $v_y$, $v_z$ are the velocities. The ratio of specific heats is given by $\gamma$ and the total energy density $\varepsilon$ is given by $\varepsilon = P/(\gamma - 1) + \rho \mathbf{v}^2/2$. Using subscripts "L" and "R" to denote the left and right states, we display the density variable for the following one-dimensional Riemann problems. All the problems in this sub-section were run on a 200 zone mesh spanning the domain $[-0.5, 0.5]$.

The first Riemann problem is the Sod problem which is given by

$(\rho_L, v_{xL}, v_{yL}, v_{zL}, P_L) = (1, 0, 0, 0, 1)$ for $x \leq 0$

$(\rho_R, v_{xR}, v_{yR}, v_{zR}, P_R) = (0.125, 0, 0, 0, 0.1)$ for $x > 0$

With $\gamma = 1.4$ the problem was run to a stopping time of 0.2 and the results are shown in Fig. 2a. We see that the contact is captured rather crisply and the other variables are free of wiggles.

The second Riemann problem is the Lax problem which is given by

$(\rho_L, v_{xL}, v_{yL}, v_{zL}, P_L) = (0.445, 0.698, 0, 0, 3.528)$ for $x \leq 0$

$(\rho_R, v_{xR}, v_{yR}, v_{zR}, P_R) = (0.5, 0, 0, 0, 0.571)$ for $x > 0$

With $\gamma = 1.4$ the problem was run to a stopping time of 0.13 and the results are shown in Fig. 2b. Again, the contact has been captured crisply.

The third test problem shows two supersonically colliding streams of fluid given by

$$(\rho_L, v_{xL}, v_{yL}, v_{zL}, P_L) = (1, 2, 0, 0, 0.2) \quad \text{for} \quad x \leq 0$$
$$(\rho_R, v_{xR}, v_{yR}, v_{zR}, P_R) = (1.5, -2, 0, 0, 0.2) \quad \text{for} \quad x > 0$$

The problem was run with $\gamma = 5/3$ to a stopping time of 0.4 and the results are shown in Fig. 2c. As before, the contact has been captured crisply and the shocks are free of wiggles.

Our last test problem shows the ability of the scheme to capture a stationary contact discontinuity. It is given by

$$(\rho_L, v_{xL}, v_{yL}, v_{zL}, P_L) = (1, 0, 0, 0, 1) \quad \text{for} \quad x \leq 0$$
$$(\rho_R, v_{xR}, v_{yR}, v_{zR}, P_R) = (0.1, 0, 0, 0, 1) \quad \text{for} \quad x > 0$$

The problem was run with $\gamma = 1.4$ to a final time of 0.25 and the results are shown in Fig. 2d. We see that the stationary contact is captured without dissipation on the computational mesh.

**VI.b) MHD Equations**

The MHD equations in one dimension are given by

$$\frac{\partial}{\partial t}\begin{pmatrix} \rho \\ \rho v_x \\ \rho v_y \\ \rho v_z \\ \varepsilon \\ B_x \\ B_y \\ B_z \end{pmatrix} + \frac{\partial}{\partial x}\begin{pmatrix} \rho v_x \\ \rho v_x^2 + P + \mathbf{B}^2/8\pi - B_x^2/4\pi \\ \rho v_x v_y - B_x B_y/4\pi \\ \rho v_x v_z - B_x B_z/4\pi \\ (\varepsilon + P + \mathbf{B}^2/8\pi) v_x - B_x (\mathbf{v} \cdot \mathbf{B})/4\pi \\ 0 \\ (v_x B_y - v_y B_x) \\ -(v_z B_x - v_x B_z) \end{pmatrix} = 0 \quad (6.2)$$

Here $\rho$ is the gas density, "P" is the gas pressure and $v_x$, $v_y$, $v_z$ are the velocities. The magnetic field components are given by $B_x$, $B_y$, $B_z$. The ratio of specific heats is given by $\gamma$ and the total energy density $\varepsilon$ is given by $\varepsilon = P/(\gamma - 1) + \rho \mathbf{v}^2/2 + \mathbf{B}^2/(8\pi)$. The MHD eigensystem has been catalogued in Roe and Balsara [46]. A very nice set of Riemann problems for MHD have been described in the papers of Brio and Wu [15], Ryu and Jones [47] and Dai and Woodward [24]; and we draw from that set of test problems. Using subscripts "*L*" and "*R*" to denote the left and right states, we display the density variable and the y-component of the magnetic field for the following one-dimensional Riemann problems. The only exception is the last MHD Riemann problem where we show the y-components of the velocity and magnetic field. All the problems in this sub-section were run on a 400 zone mesh spanning the domain $[-0.5, 0.5]$.

Our first test problem is from Brio and Wu and is given by

$$\left(\rho_L, v_{xL}, v_{yL}, v_{zL}, P_L, B_{yL}, B_{zL}\right) = (1, 0, 0, 0, 1, \sqrt{4\pi}, 0) \quad \text{for } x \leq 0$$

$$\left(\rho_R, v_{xR}, v_{yR}, v_{zR}, P_R, B_{yR}, B_{zR}\right) = (0.125, 0, 0, 0, 0.1, -\sqrt{4\pi}, 0) \quad \text{for } x > 0$$

We set $\gamma = 2$ and $B_x = 0.75\sqrt{4\pi}$, and we ran this problem to a final time of 0.1. The results are shown in Figs. 3a and 3b. We see that the compound wave is properly captured and the density profile in the contact discontinuity is crisp.

Our second test problem is from Ryu and Jones and shows the formation of all seven waves in an MHD Riemann problem. It is given by

$$\left(\rho_L, v_{xL}, v_{yL}, v_{zL}, P_L, B_{yL}, B_{zL}\right) = (1.08, 1.2, 0.01, 0.5, 0.95, 3.6, 2.0) \quad \text{for } x \leq 0$$

$$\left(\rho_R, v_{xR}, v_{yR}, v_{zR}, P_R, B_{yR}, B_{zR}\right) = (1, 0, 0, 0, 1, 4, 2) \quad \text{for } x > 0$$

We used $\gamma = 5/3$ and $B_x = 2$, and we ran this test problem to a final time of 0.2. The results are shown in Figs. 3c and 3d. The fast and slow magnetosonic shocks reveal themselves in the density profile. The density also shows the contact discontinuity in between the two slow magnetosonic shocks. The two Alfven waves reveal themselves more clearly in the y-component of the magnetic field. Since the entropy wave and Alfven waves are linearly degenerate, we see that they have been captured very crisply by the HLLI-GRP solver.

Our third test problem is from Dai and Woodward and shows the effect of two MHD streams that collide supersonically with each other. It is given by

$$\left(\rho_L, v_{xL}, v_{yL}, v_{zL}, P_L, B_{yL}, B_{zL}\right) = (0.15, 21.55, 1, 1, 0.28, -2, -1) \quad \text{for } x \leq 0$$

$$\left(\rho_R, v_{xR}, v_{yR}, v_{zR}, P_R, B_{yR}, B_{zR}\right) = (0.1, -26.45, 0, 0, 0.1, 2, 1) \quad \text{for } x > 0$$

We ran this problem with $\gamma = 5/3$ and $B_x = 0$ to a final time of 0.04. The results are shown in Figs. 3e and 3f. We see the formation of two extremely high Mach number fast shocks. Because the longitudinal magnetic field is zero, the slow shocks do not form.

Our fourth test problem is from Ryu and Jones and shows the formation of a switch-on fast magnetosonic shock. It is given by

$$\left(\rho_L, v_{xL}, v_{yL}, v_{zL}, P_L, B_{yL}, B_{zL}\right) = (1, 0, 0, 0, 1, \sqrt{4\pi}, 0) \quad \text{for } x \leq 0$$

$$\left(\rho_R, v_{xR}, v_{yR}, v_{zR}, P_R, B_{yR}, B_{zR}\right) = (0.2, 0, 0, 0, 0.1, 0, 0) \quad \text{for } x > 0$$

We ran this problem with $\gamma = 5/3$ and $B_x = \sqrt{4\pi}$ to a final time of 0.15. We see from Figs. 3g and 3h that the transverse component of the magnetic field is indeed zero before it gets run over by the right-going fast shock and its value increases in the post-shock region. For this reason, the right-going fast shock is thought to "switch on" the transverse component of the magnetic field. We see that the switch-on shock is properly captured by our method.

Our fifth test problem shows the ability of the method to retain a stationary Alfven wave on the computational mesh. It is given by

$$\left(\rho_L, v_{xL}, v_{yL}, v_{zL}, P_L, B_{yL}, B_{zL}\right) = (1/(4\pi), -1, 1, -1, 1, -1, 1) \quad \text{for } x \leq 0$$

$$\left(\rho_R, v_{xR}, v_{yR}, v_{zR}, P_R, B_{yR}, B_{zR}\right) = (1/(4\pi), -1, -1, -1, 1, 1, 1) \quad \text{for } x > 0$$

We ran this problem with $\gamma = 1.4$ and $B_x = 1$ to a final time of 0.1. Figs. 3i and 3j show the y-velocity and the y-magnetic field. We see that the Alfven wave is very crisply captured on the mesh. This highlights the ability of our HLLI-GRP method to capture any linearly degenerate characteristic field with minimal dissipation.

### VI.c) Shallow Water Equations with Non-Constant Bathymetry

Let us consider the augmented single-layer shallow water equations with time-independent bottom topography given by

$$\frac{\partial h}{\partial t} + \frac{\partial (hu)}{\partial x} = 0,$$

$$\frac{\partial (hu)}{\partial t} + \frac{\partial}{\partial x}\left(hu^2 + \frac{1}{2}gh^2\right) + gh\frac{\partial b}{\partial x} = S_{fX},$$

$$\frac{\partial (hv)}{\partial t} + \frac{\partial}{\partial x}(huv) = S_{fY},$$

$$\frac{\partial b}{\partial t} = 0,$$

(6.3)

where $h$ is the water depth, $u$ is the normal velocity, $v$ is the transverse velocity, $b = b(x)$ is the bottom topography, $g$ is the gravity acceleration, $S_{fX}$ and $S_{fY}$ are source terms which represent the bottom friction components. In this paper we are interested in the classical Manning formulation (Manning [43] and Chertock *et al.* [22])

$$S_{fX} = -gn^2 \frac{u\sqrt{u^2+v^2}}{h^{1/3}}; \quad S_{fY} = -gn^2 \frac{v\sqrt{u^2+v^2}}{h^{1/3}}, \tag{6.4}$$

where $n$ is the Manning coefficient. Notice that if $h$ takes small values these source terms became stiff. In a matrix form the system (6.3) takes the form of the system in eqn. (5.1) where

$$\mathbf{U} = \begin{pmatrix} h \\ hu \\ hv \\ b \end{pmatrix} \; ; \; \mathbf{F}(\mathbf{U}) = \begin{pmatrix} hu \\ hu^2 + \frac{1}{2}gh^2 \\ huv \\ 0 \end{pmatrix} \; ; \; \mathbf{B}(\mathbf{U}) = \begin{pmatrix} 0 & 0 & 0 & 0 \\ 0 & 0 & 0 & gh \\ 0 & 0 & 0 & 0 \\ 0 & 0 & 0 & 0 \end{pmatrix} \; ; \; \mathbf{S}(\mathbf{U}) = \begin{pmatrix} 0 \\ S_{fX} \\ S_{fY} \\ 0 \end{pmatrix} \tag{6.5}$$

Notice that in the case $n = 0$, the system takes the form of eqn. (3.1). Furthermore, the system of equations (6.3) correspond to a hyperbolic system with characteristic matrix given by

$$\mathbf{A}(\mathbf{U}) = \begin{pmatrix} 0 & 1 & 0 & 0 \\ gh - u^2 & 2u & 0 & gh \\ -uv & v & u & 0 \\ 0 & 0 & 0 & 0 \end{pmatrix} \tag{6.6}$$

and eigenvalues; $\lambda_1 = u - c$, $\lambda_2 = 0$, $\lambda_3 = u$, $\lambda_4 = u + c$, where $c = \sqrt{gh}$. In the following paragraphs we provide the set of orthonormalized left and right eigenvectors for this system.

The eigenvalue $\lambda_1$ is the left-going wave which corresponds to a genuinely nonlinear field and its orthonormalized right and left eigenvectors are given by

$$r_1^p = (1, \; u-c, \; v, \; 0)^T$$
$$l_1^p = \left( \frac{c+u}{2c}, \; -\frac{1}{2c}, \; 0, \; -\frac{c}{2(u-c)} \right) \tag{6.7}$$

The superscript "$p$" denotes that the eigenvectors are specified in terms of primitive variables.

The eigenvalue $\lambda_4$ is the right-going wave which corresponds to a genuinely nonlinear field and its orthonormalized right and left eigenvectors are given by

$$r_4^p = \begin{pmatrix} 1, & u+c, & v, & 0 \end{pmatrix}^T$$
$$l_4^p = \begin{pmatrix} \dfrac{c-u}{2c}, & \dfrac{1}{2c}, & 0, & \dfrac{c}{2(u+c)} \end{pmatrix} \tag{6.8}$$

The eigenvalue $\lambda_2$ corresponds to a stationary wave associated with the bottom jump and this corresponds to a linearly degenerate field. The orthonormalized right and left eigenvectors are given by

$$r_2^p = \begin{pmatrix} 1, & 0, & v, & \dfrac{(u^2-c^2)}{c^2} \end{pmatrix}^T$$
$$l_2^p = \begin{pmatrix} 0, & 0, & 0, & \dfrac{c^2}{(u^2+c^2)} \end{pmatrix} \tag{6.9}$$

The eigenvalue $\lambda_3$ corresponds to a shear wave associated with the transverse flow velocity $v$ and this corresponds to a linearly degenerate field. The orthonormalized right and left eigenvectors are given by

$$r_3^p = \begin{pmatrix} 0, & 0, & 1, & 0 \end{pmatrix}^T$$
$$l_3^p = \begin{pmatrix} -v, & 0, & 1, & 0 \end{pmatrix} \tag{6.10}$$

For the sake of completeness, we also document that the Jacobian of the source term has the form

$$\dfrac{\partial \mathbf{S}}{\partial \mathbf{U}}(\mathbf{U}) = \begin{pmatrix} 0 & 0 & 0 & 0 \\ \dfrac{gn^2 u\sqrt{u^2+v^2}}{h^{1/3}} & -\dfrac{gn^2(2u^2+v^2)}{h^{1/3}\sqrt{u^2+v^2}} & -\dfrac{gn^2 uv}{h^{1/3}\sqrt{u^2+v^2}} & 0 \\ \dfrac{gn^2 v\sqrt{u^2+v^2}}{h^{1/3}} & -\dfrac{gn^2 uv}{h^{1/3}\sqrt{u^2+v^2}} & -\dfrac{gn^2(u^2+2v^2)}{h^{1/3}\sqrt{u^2+v^2}} & 0 \\ 0 & 0 & 0 & 0 \end{pmatrix} \tag{6.11}$$

The first set of tests corresponds to Riemann problems reported in Dumbser and Balsara [28]; the initial conditions are defined in Table 1, All the Riemann problems were run on a 100 zone mesh spanning domains in the form $[x_L, x_R]$ and we have set the CFL number to 0.9. Figure 4 shows the results. Riemann problem RP0 consists on a flat free surface with a jump discontinuity in the bottom topography and a shear wave. These profiles are stationary solutions, which are well captured by the HLLI-GRP scheme. Riemann problems PR1 and RP2 correspond to dam-break problems. Riemann problem RP1 faces a dry bed case, which is solved very well. The remaining Riemann problems are also properly solved in agreement with the results reported in Dumbser and Balsara [28].

**Table 1 Single-layer shallow water equations; left and right initial states ( water depth $h$, velocity components $u$ and $v$, bottom height $b$ ), final output time, computational domain $[x_L, x_R]$ and initial position $x_C$ of the discontinuity.**

| Case | $h_L$ | $u_L$ | $v_L$ | $b_L$ | $h_R$ | $u_R$ | $v_R$ | $b_R$ | $t_{end}$ | $x_L$ | $x_R$ | $x_C$ |
|------|-------|-------|-------|-------|-------|-------|-------|-------|-----------|-------|-------|-------|
| RP0 | 2 | 0 | 1 | 0 | 1 | 0 | −1 | 1 | 1 | 0 | 1 | 0.5 |
| RP1 | 1 | 0 | 0 | 0 | $10^{-14}$ | 0 | 0 | 0 | 0.075 | 0 | 1 | 0.5 |
| RP2 | 1.46184 | 0 | 0 | 0 | 0.30873 | 0 | 0 | 0.2 | 1 | −5 | 5 | 0 |
| RP3 | 0.75 | −9.49365 | 0 | 0 | 1.10594 | −4.94074 | 0 | 0.2 | 1 | −15 | 5 | 0 |
| RP4 | 0.75 | −1.35624 | 0 | 0 | 1.10594 | −4.94074 | 0 | 0.2 | 1 | −10 | 4 | 0 |

The second set of problems corresponds to small perturbation of steady flow over a slanted surface, Table 2 shows the parameters. The initial condition for all these problems has the form

$$h(x,0) = h_0 + \begin{cases} 0.2h_0, & 1 \leq x \leq 1.25, \\ 0, & \text{otherwise} \end{cases} \qquad q(x,0) = u(x,0)h(x,0) = q_0$$

It can be seen that the initial condition is a small perturbation of $h(x, 0) = h_0$. All tests were run on a 100 zone mesh spanning the domain [0,25] and we have set the CFL number to 0.9. Figure 4 shows the results. The perturbation travels to the right leaving the domain for large times. The initial profile modifies its shape as it advances and it reaches a steady state after a large period of time. This behavior is in agreement with that reported in the literature by Chertock *et al*. [22].

**Table 2 Single-layer shallow water equations with friction source term; set of data for small perturbations of steady flow over a slanted surface (the water depth $h_0$, the water discharge $q_0$, constant slope of bottom topography $b_x$, the Manning coefficient $n$).**

| Test | $h_0$ | $q_0$ | $n$ | $b_x$ |
|------|-------|-------|-----|-------|
| 1 | 0.09564 | 0.1 | 0.02 | -0.01 |
| 2 | 0.02402 | 0.002 | 0.1 | -0.01 |
| 3 | 0.44894 | 2 | 0.1 | $-1/\sqrt{3}$ |

### VI.d) Compressible Navier Stokes Equations in Relaxation form

The one-dimensional compressible Navier-Stokes equations are given by

$$\frac{\partial \rho}{\partial t} + \frac{\partial (\rho u)}{\partial x} = 0,$$

$$\frac{\partial (\rho u)}{\partial t} + \frac{\partial}{\partial x}\left(\rho u^2 + p\right) = \frac{4}{3}\frac{\partial}{\partial x}\left(\mu \frac{\partial}{\partial x}\right), \quad (6.12)$$

$$\frac{\partial E}{\partial t} + \frac{\partial}{\partial x}\left(u(E+p)\right) = \frac{4}{3}\frac{\partial}{\partial x}\left(\mu u \frac{\partial}{\partial x}\right) - \frac{\partial}{\partial x}(\kappa T),$$

where $\rho$ is the density, $u$ is the velocity, $E$ is the total energy, $p$ is the pressure, $T$ is the temperature, $\mu$ is the viscosity coefficient, $\kappa$ is the heat transfer coefficient. The total energy is given by

$$E = \rho\left(\frac{1}{2}u^2 + e(p,\rho)\right), \quad (6.13)$$

where $e(p,\rho)$ is the specific internal energy, which in this work is chosen as

$$e(p,\rho) = \frac{p}{(\gamma-1)\rho}, \quad (6.14)$$

with $\gamma$ the ratio of specific heats. The temperature is given by $T = \frac{p}{R\rho}$, where $R$ is the gas constant. For the viscosity we use Sutherland's law given by

$$\mu(T) = \mu_0 \cdot \left(\frac{T}{T_0}\right)^\beta \cdot \left(\frac{T+s}{T_0+s}\right) \quad (6.15)$$

Where $\mu_0$, $T_0$, $\beta$ and $s$ are constant. The heat conduction is related to the viscosity as

$\kappa = \frac{\mu \gamma c_v}{\Pr}$, here $\Pr$ is the Prandtl number and $c_v = \frac{R}{(\gamma-1)}$ is the capacity at constant volume. This completes the description of the one-dimensional compressible Navier-Stokes.

Now, for solving this system numerically, we reformulate it as a hyperbolic system by using a Cattaneo [20], [21] type relaxation proposed by Toro and Montecinos [44]. So we introduce two auxiliary variables $\psi_1$ and $\psi_2$ with the evolutionary equations

$$\frac{\partial \psi_1}{\partial t} = \frac{1}{\varepsilon}\left(\frac{\partial u}{\partial x} - \psi_1\right), \quad \frac{\partial \psi_2}{\partial t} = \frac{1}{\varepsilon}\left(\frac{\partial T}{\partial x} - \psi_2\right). \quad (6.16)$$

Notice that $\psi_1 \to \frac{\partial u}{\partial x}$, $\psi_2 \to \frac{\partial T}{\partial x}$, as $\varepsilon \to 0$. In this form we build a new system, which has the form of the system (3.1) but in a conservation form, ($\mathbf{B} = 0$), the vector of state, the flux and source are given by

$$\mathbf{U} = \begin{pmatrix} \rho \\ \rho u \\ E \\ \psi_1 \\ \psi_2 \end{pmatrix}, \quad \mathbf{F}(\mathbf{U}) = \begin{pmatrix} \rho u \\ \rho u^2 + p - \frac{4\mu\psi_1}{3} \\ u\left(E + p - \frac{4\mu\psi_1}{3}\right) + \kappa\psi_2 \\ -u/\varepsilon \\ -T/\varepsilon \end{pmatrix}, \quad \mathbf{S}(\mathbf{U}) = \begin{pmatrix} 0 \\ 0 \\ 0 \\ -\psi_1/\varepsilon \\ -\psi_2/\varepsilon \end{pmatrix} \tag{6.17}$$

For this system, the eigenstructure cannot be obtained explicitly; only partial information can be obtained for this system. In Montecinos and Toro (2014) the authors have identified a necessary condition for the case $\kappa = 0$ which guarantees that the relaxation system has real eigenvalues, it is the case of interest in this section. The eigenvectors are approximated using numerical recipes for obtaining the eigenvectors.

The test problems are those reported in Montecinos and Toro [44], the initial condition is given by

$$(\rho_L, u_L, p_L) = (1.29, 0, 2929.73), \quad \text{for } x \leq 0$$

$$(\rho_R, u_R, p_R) = (1.784, 0, 4349.31), \quad \text{for } x \geq 0$$

All problems in this section were run on a 100 zone mesh spanning the domain $[-1, 1]$ and stopped at the output time $t = 0.01$. We have set the CFL number to 0.7, the relaxation parameter has been chosen as $\varepsilon = 10^{-4}$ and $\gamma = 1.4$.

Figures 6, 7 and 8 show the results for $\mu = 2$, $\mu = 0.2$. and $\mu = 0.01$, respectively. The velocity profile evidences the influence of the viscosity, thus the reformulation obtained from the Cattaneo relaxation procedure is able to reproduce this behavior. For the chosen value of $\varepsilon$ the system is hyperbolic, furthermore, the system is stiff. The HLLI-GRP is able to provide very good agreement with respect to the reported solutions in Montecinos and Toro [44].

## VII) Conclusions

While the generalized Riemann problem has been extensively studied for certain select systems of hyperbolic equations, a general purpose and easily-implementable strategy for obtaining the GRP for any hyperbolic system has been missing in the research literature. This deficiency stems from the fact that the underlying Riemann problems can be vastly different for different hyperbolic systems. If a GRP solver is to be built out of an exact Riemann problem solver, we will necessarily be forced to develop GRP solvers for each hyperbolic system on a case-by-case basis.

The recently-proposed HLLI Riemann solver is an approximate Riemann solver that is indeed universal, in that it is applicable to any hyperbolic system, whether in conservation form or with non-conservative products. The HLLI Riemann solver is also complete in the sense that if it is given a complete set of eigenvectors, it represents all waves with minimal dissipation. It is, therefore, very desirable to upgrade it so as to obtain an HLLI-GRP solver. We have undertaken this task in the present paper. Consequently, we present an HLLI-GRP solver that inherits all the favorable properties of the HLLI parent. Section IV.a provides all the implementation-related details for implementing an HLLI-GRP solver for hyperbolic systems in conservation form. Section IV.b provides similar implementation-related details for implementing the HLLI-GRP solver for hyperbolic systems with non-conservative products.

A Riemann solver does not necessarily need to be changed when stiff source terms are present. However, a GRP solver needs to be modified to accommodate the presence of stiff source terms. In this paper we also show how our GRP solver can be adapted to the solution of hyperbolic systems with stiff source terms. Section V.b provides all the implementation-related details for implementing the HLLI-GRP solver for hyperbolic systems in conservation form with stiff source terms. Section V.c provides all the implementation-related details for implementing the HLLI-GRP solver for hyperbolic systems with non-conservative products including stiff source terms.

Results from several stringent test problems are shown. These test problems are drawn from many different hyperbolic systems, and include hyperbolic systems in conservation form; with non-conservative products; and with stiff source terms. The present generalized Riemann problem solver performs well on all of them.

**Acknowledgements**

DSB thanks the hospitality of Institute of Applied Physics and Computational Mathematics, Beijing when this work started. Jiequan Li is supported by NSFC with nos. 11771054 and 11371063, and by Foundation of LCP at Institute of Applied Physics and Computational Mathematics, Beijing.

DSB also acknowledges support via NSF grants NSF-ACI-1533850, NSF-DMS-1622457, NSF_ACI-1713765. Several simulations were performed on a cluster at UND that is run by the Center for Research Computing. Computer support on NSF's XSEDE and Blue Waters computing resources is also acknowledged.

**Figure Captions**

*Fig. 1 corresponds to eqn. (3.4) and shows how the states in the HLL Riemann solver can be parameterized in terms of the similarity variable. This parametrization is very useful for hyperbolic systems with non-conservative products.*

*Fig. 2 shows the density profile from four Riemann problems involving Euler flow. Figs. 2a and 2b show the results from the Sod and Lax problems. Fig. 2c shows the result from a test problem with supersonically colliding fluids. Fig. 2d shows that an isolated contact is preserved exactly on the mesh.*

*Fig. 3 shows the result of several MHD Riemann problems. In all examples, except the last one, we show the density and the y-component of the magnetic field. Figs. 3a and 3b show the results from the Brio-Wu test problem. Figs. 3c and 3d show the results from a seven-wave test problem from Ryu and Jones. Figs. 3e and 3f show the results of a problem with colliding MHD streams from Dai and Woodward. Figs. 3g and 3h show the results of a Riemann problem with a switch-on fast shock. Figs. 3i and 3j show the y-velocity and y-magnetic field of a standing Alfven discontinuity, showing that it is captured without dissipation on the mesh.*

*Fig. 4 shows the result of the shallow water Riemann problems. Figs. 4a and 4b show the density and velocity for the Riemann problem RP0. Figs. 4c and 4d show the free surface for the Riemann problems RP1 and RP2, respectively. Figs. 4e and 4f show the free surface for the Riemann problems RP3 and RP4, respectively.*

*Fig. 5 shows the free surface and the topography for the shallow water tests for flow over a slanted surface. Figs. 5a, 5b and 5c show the profiles for the test 1 at time $t = 0, t = 1$ and $t = 100$, respectively. Figs. 5d, 5e and 5f show the profiles for the test 2 at time $t = 0, t = 0.5$ and $t = 100$, respectively. Figs. 5g, 5h and 5i show the profiles for the test 3 at time $t = 0, t = 0.05$ and $t = 100$, respectively.*

*Fig. 6 shows the results for the compressible Navier-Stokes flow with $\mu = 2Pa/s$ at the output time 0.01. Fig. 6a shows the profile for the density. Fig. 6b shows the profile for the x-component of the velocity. Fig. 6c shows the profile for the density. Fig. 6d shows the internal energy.*

*Fig. 7 shows the results for the compressible Navier-Stokes flow with $\mu = 0.2Pa/s$ at the output time 0.01. Fig. 7a shows the profile for the density. Fig. 7b shows the profile for the x-component of the velocity. Fig. 7c shows the profile for the density. Fig. 7d shows the internal energy.*

*Fig. 8 shows the results for the compressible Navier-Stokes flow with $\mu = 0.001Pa/s$ at the output time 0.01. Fig. 8a shows the profile for the density. Fig. 8b shows the profile for the x-component of the velocity. Fig. 8c shows the profile for the density. Fig. 8d shows the internal energy.*

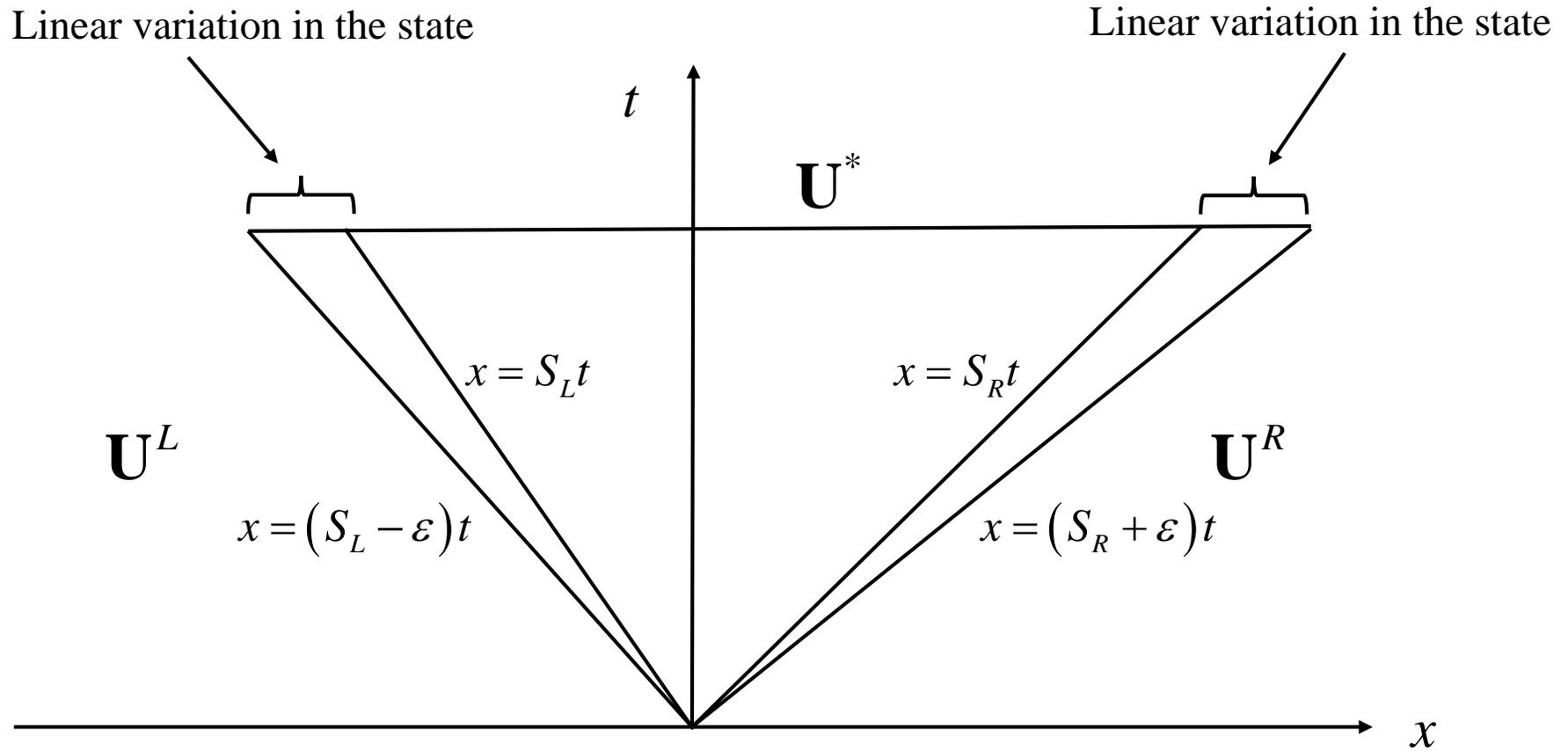

Fig. 1 corresponds to eqn. (3.4) and shows how the states in the HLL Riemann solver can be parameterized in terms of the similarity variable. This parametrization is very useful for hyperbolic systems with non-conservative products.

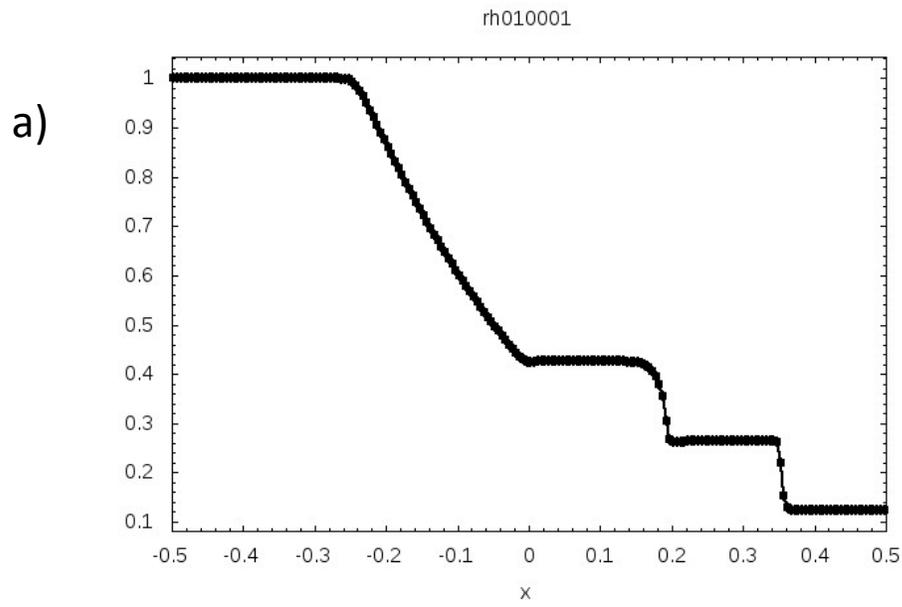 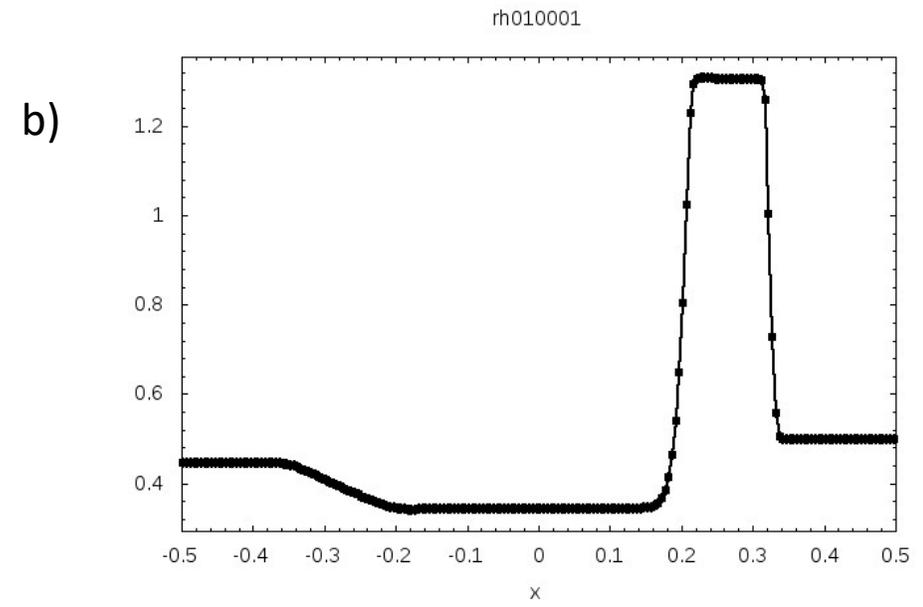
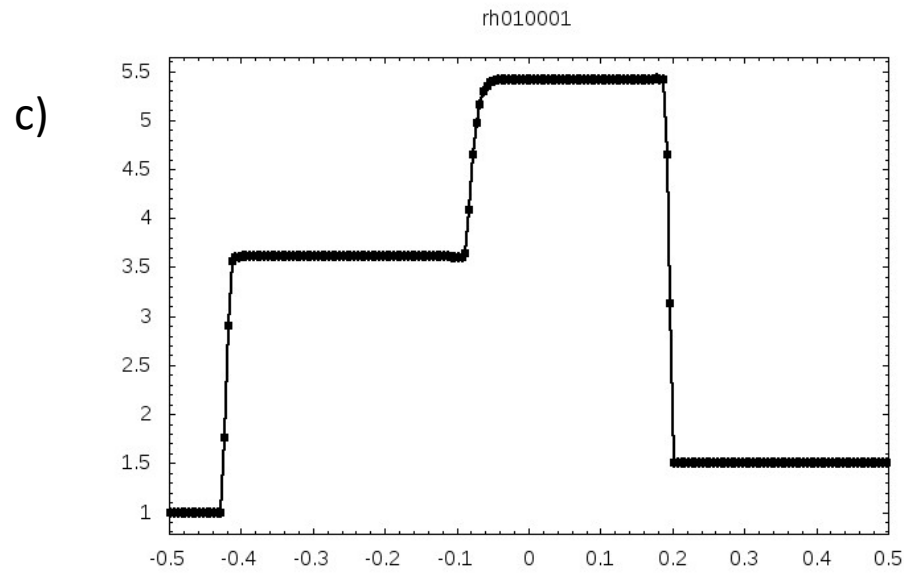 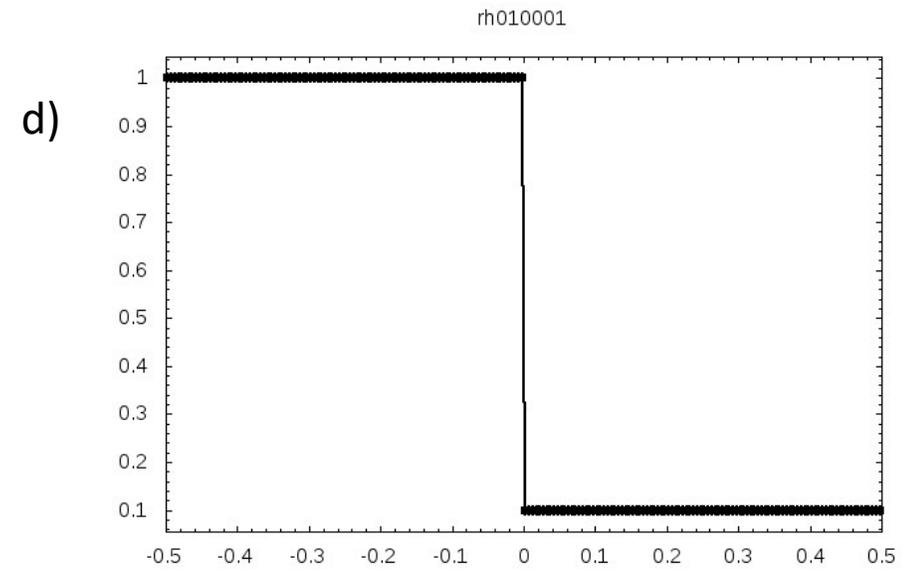

*Fig. 2 shows the density profile from four Riemann problems involving Euler flow. Figs. 2a and 2b show the results from the Sod and Lax problems. Fig. 2c shows the result from a test problem with supersonically colliding fluids. Fig. 2d shows that an isolated contact is preserved exactly on the mesh.*

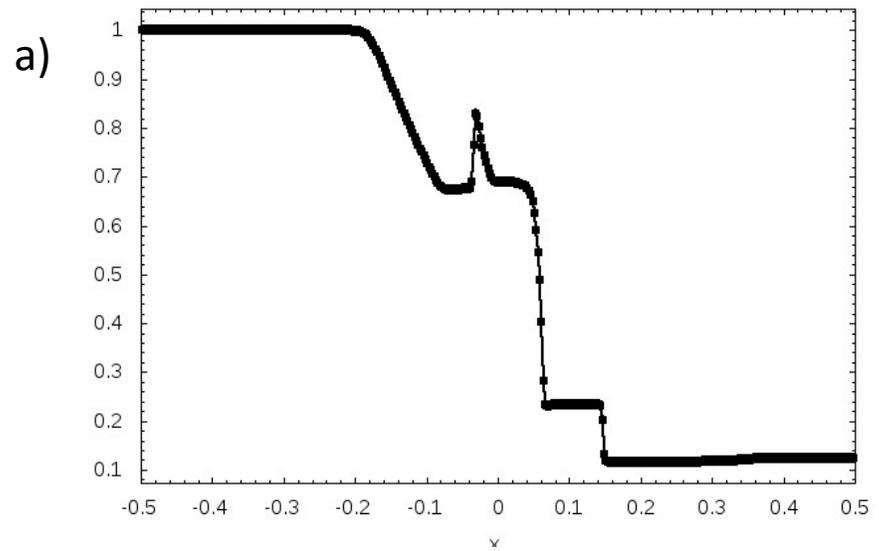
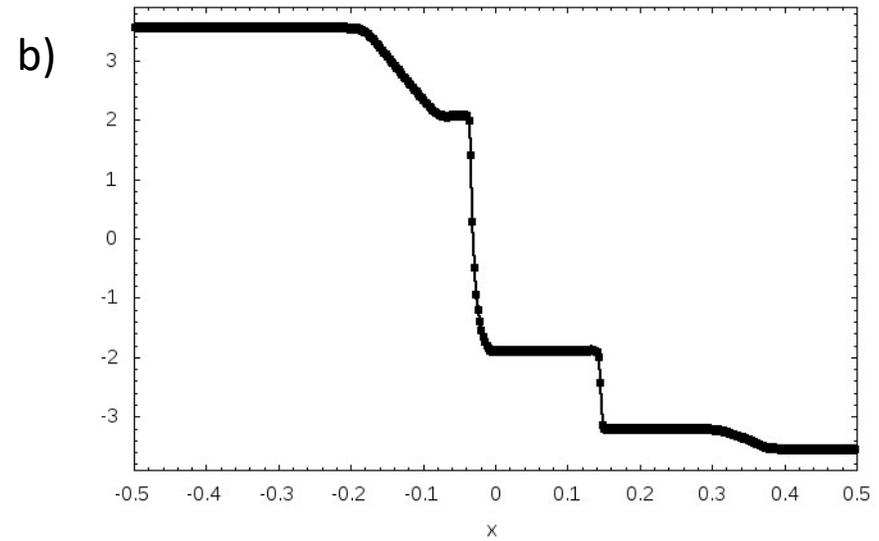
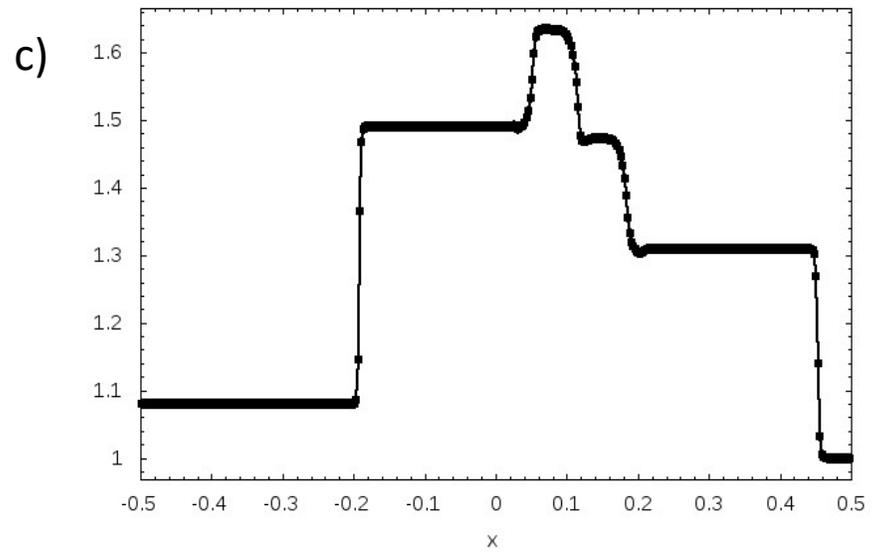
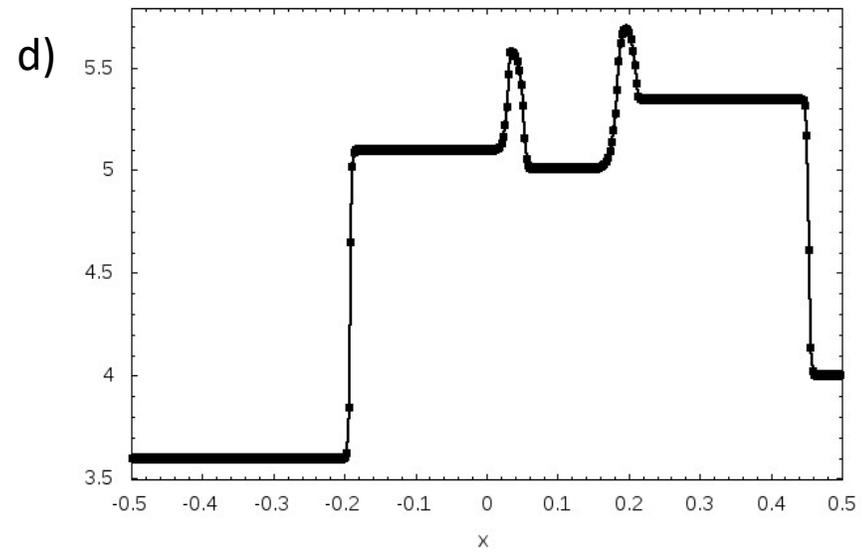

e) 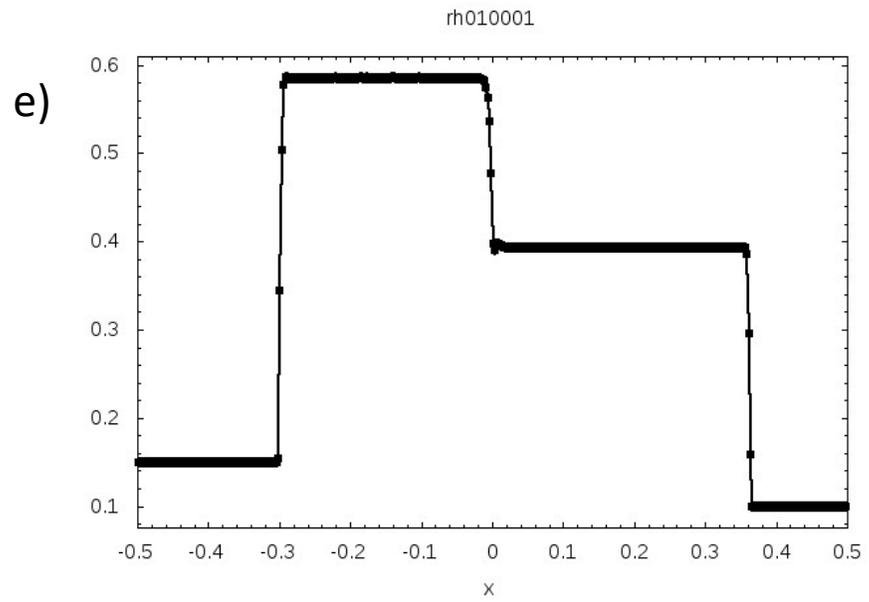

f) 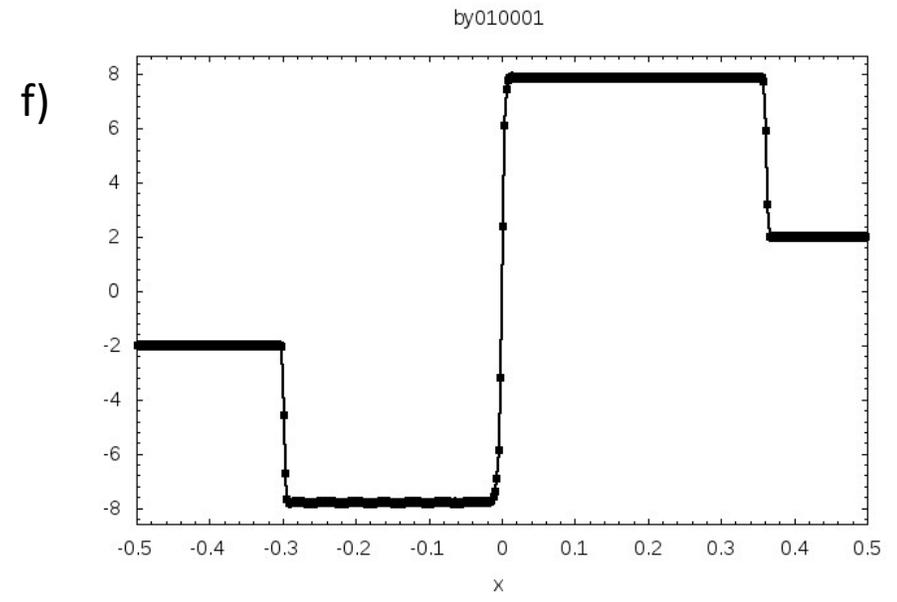

g) 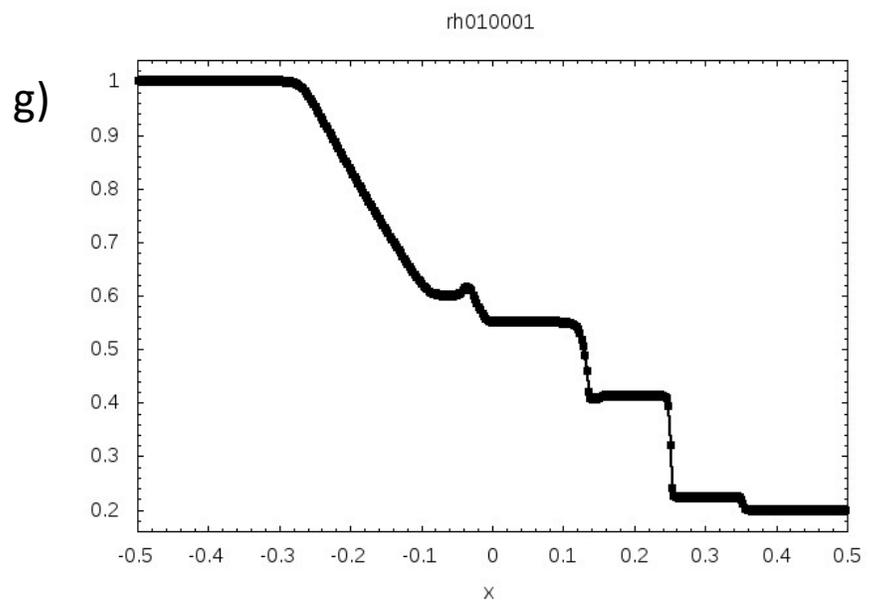

h) 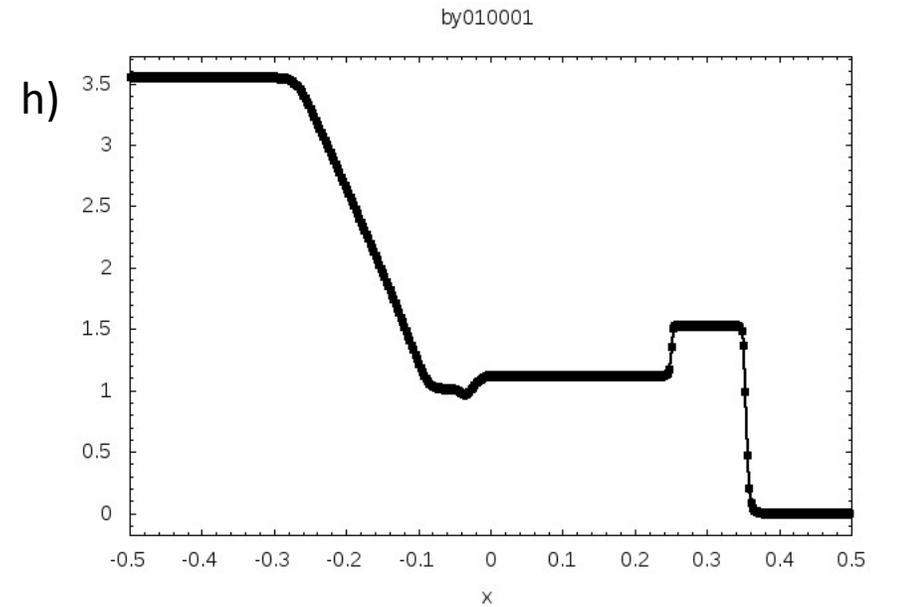

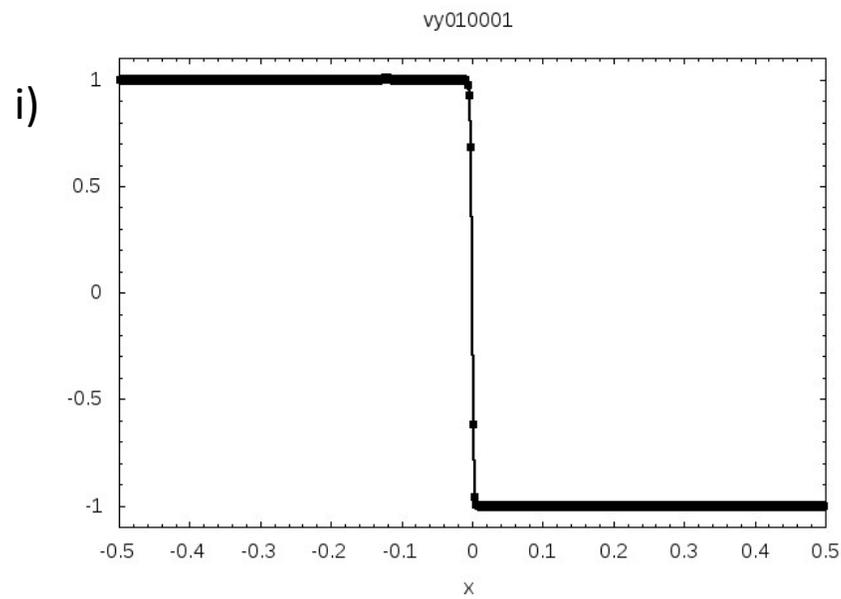 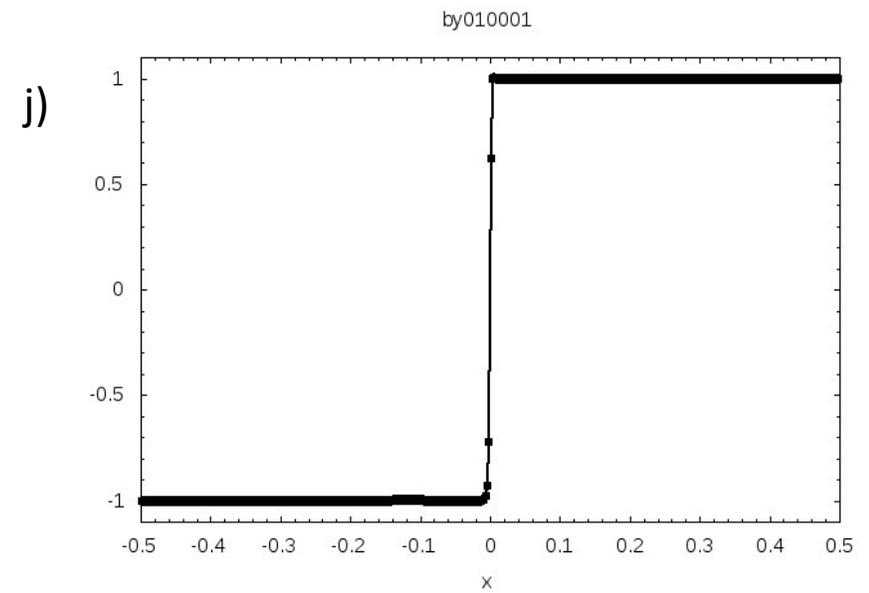

Fig. 3 shows the result of several MHD Riemann problems. In all examples, except the last one, we show the density and the y-component of the magnetic field. Figs. 3a and 3b show the results from the Brio-Wu test problem. Figs. 3c and 3d show the results from a seven-wave test problem from Ryu and Jones. Figs. 3e and 3f show the results of a problem with colliding MHD streams from Dai and Woodward. Figs. 3g and 3h show the results of a Riemann problem with a switch-on fast shock. Figs. 3i and 3j show the y-velocity and y-magnetic field of a standing Alfven discontinuity, showing that it is captured without dissipation on the mesh.

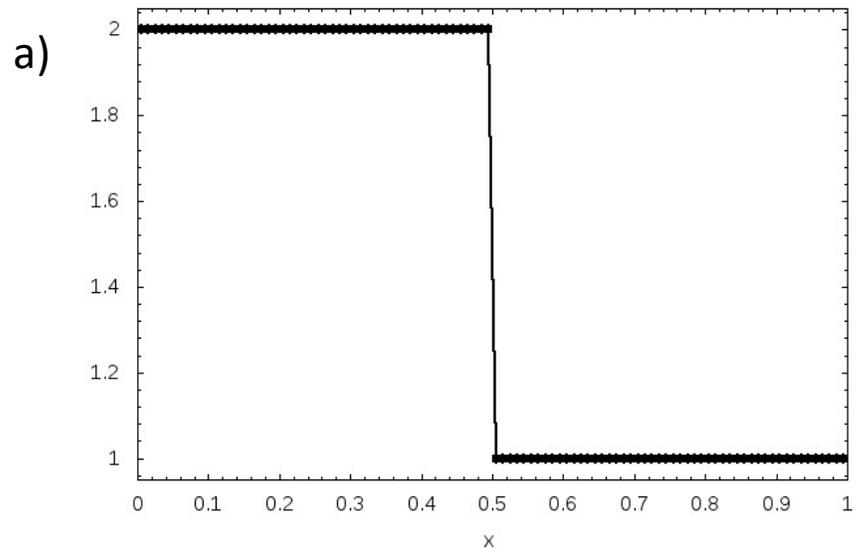
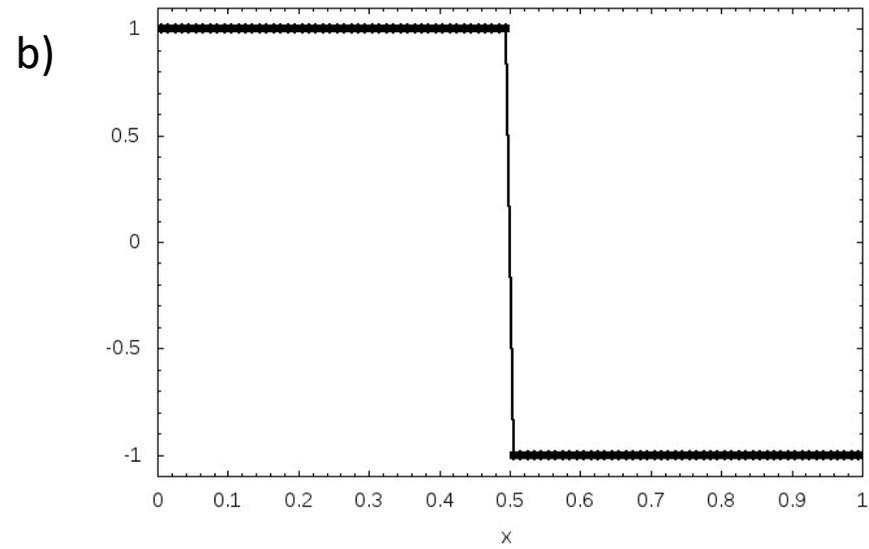
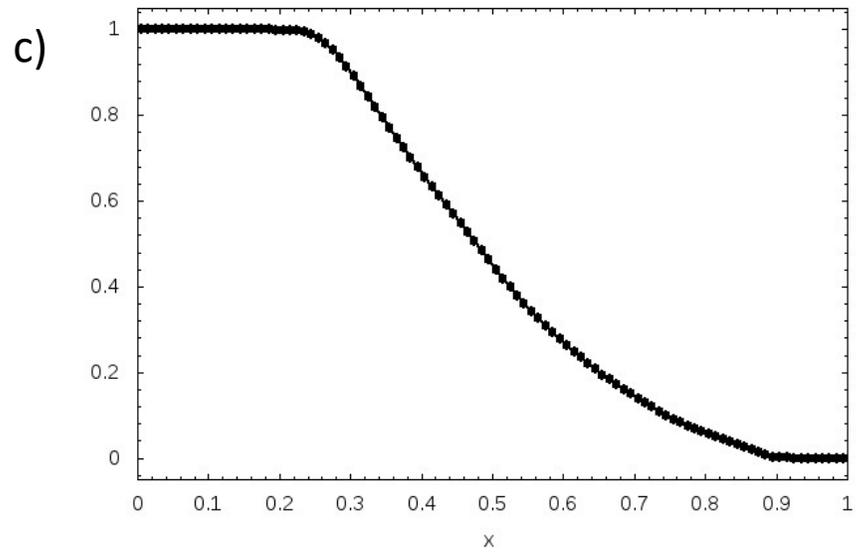
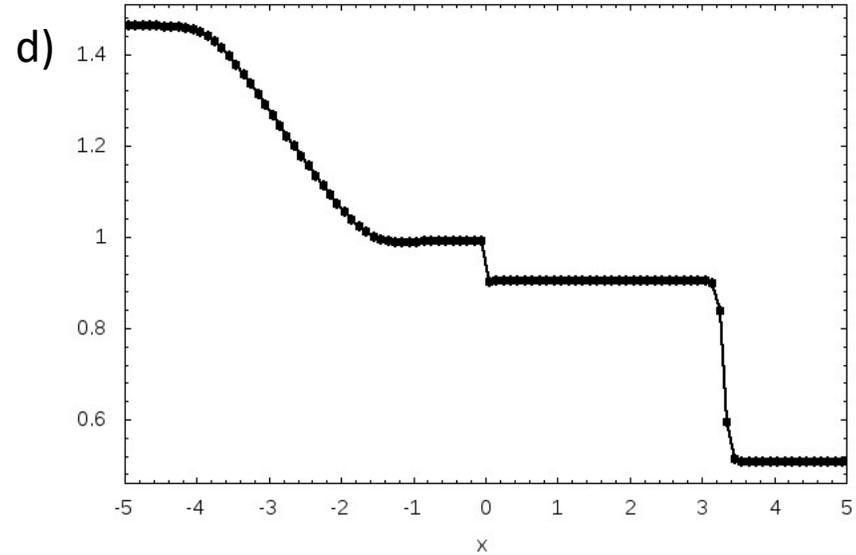

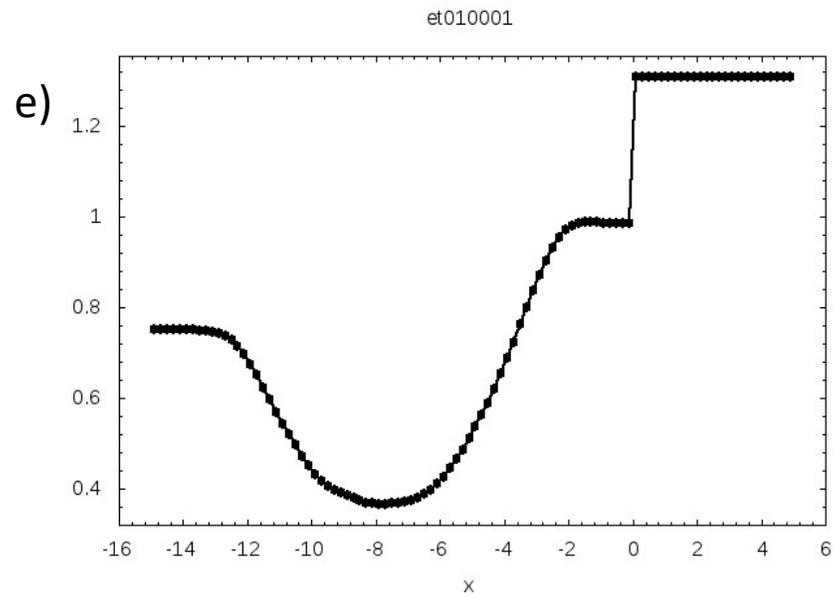 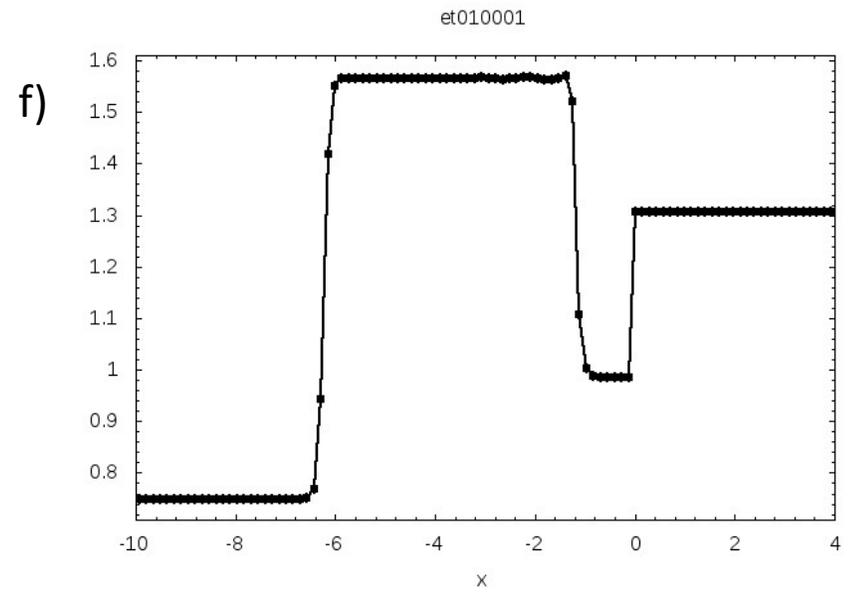

Fig. 4 shows the result of the shallow water Riemann problems. Figs. 4a and 4b show the density and velocity for the Riemann problem RP0. Figs. 4c and 4d show the free surface for the Riemann problems RP1 and RP2, respectively. Figs. 4e and 4f show the free surface for the Riemann problems RP3 and RP4, respectively.

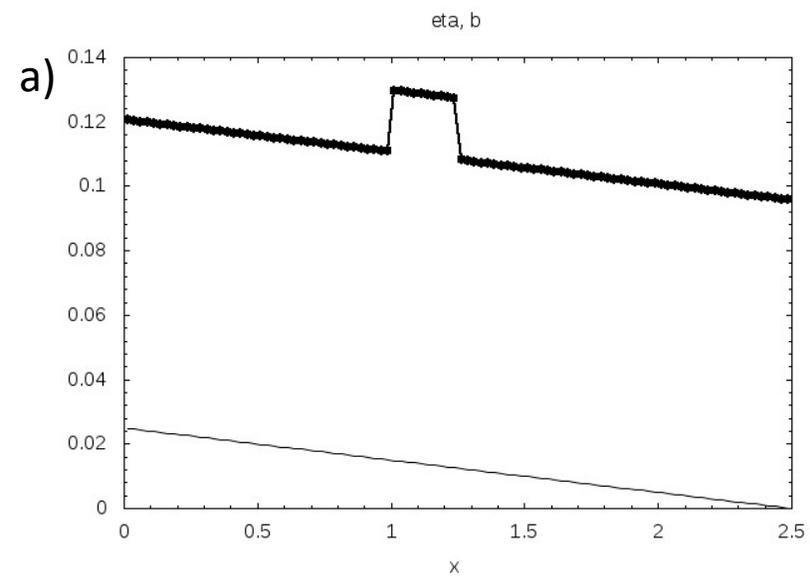 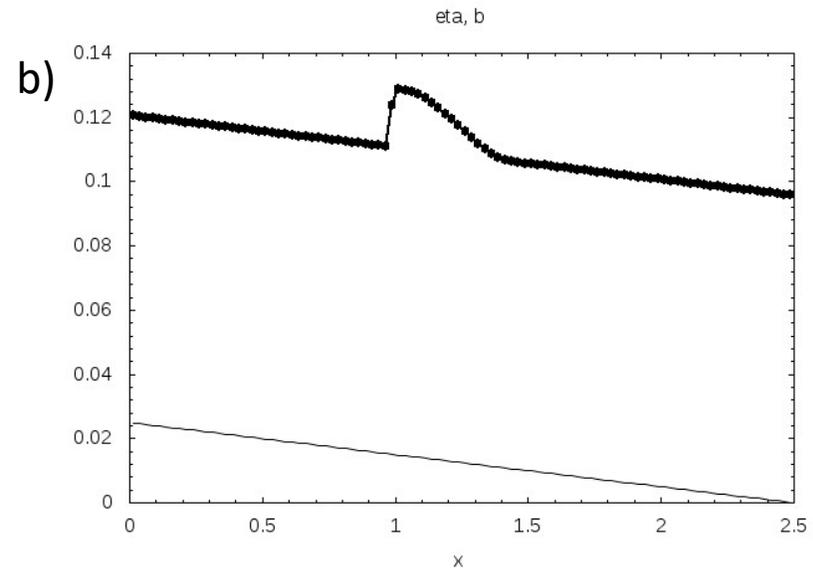 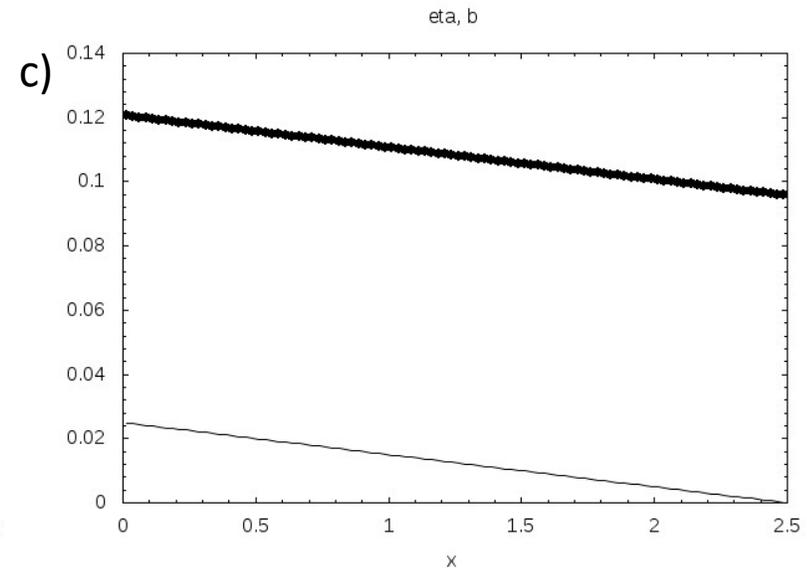
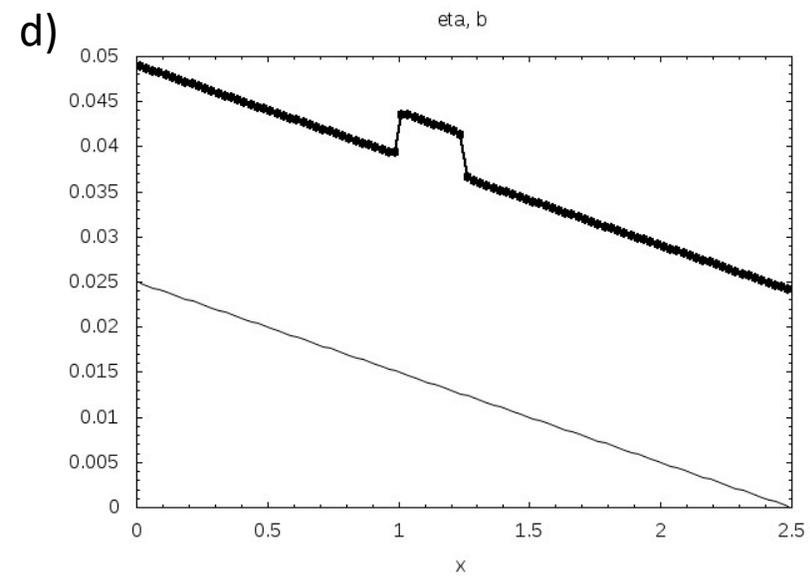 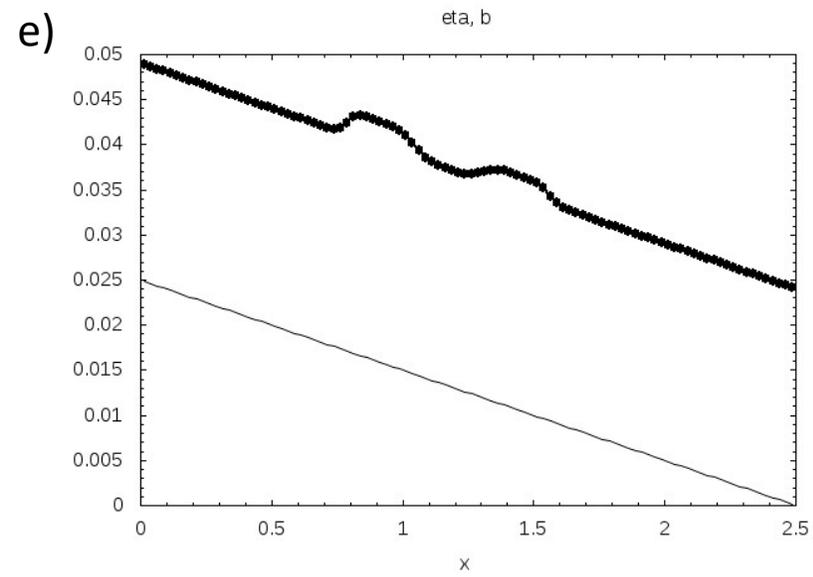 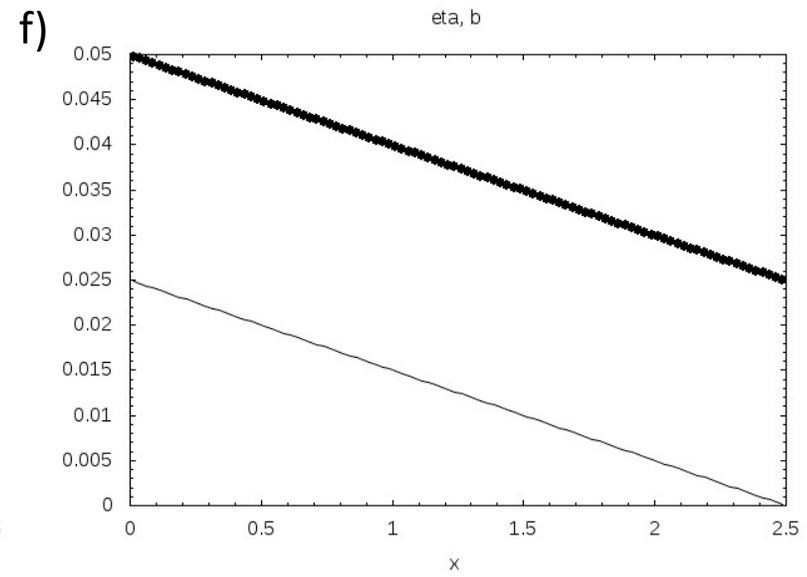

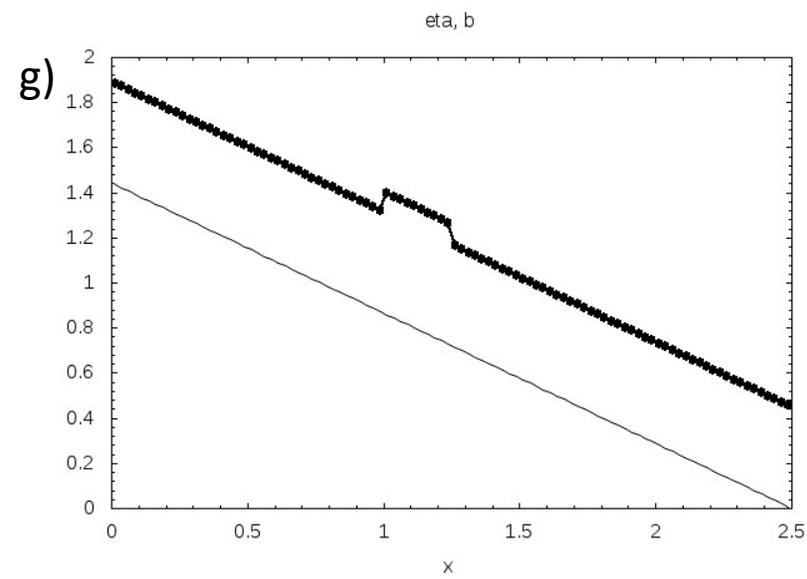 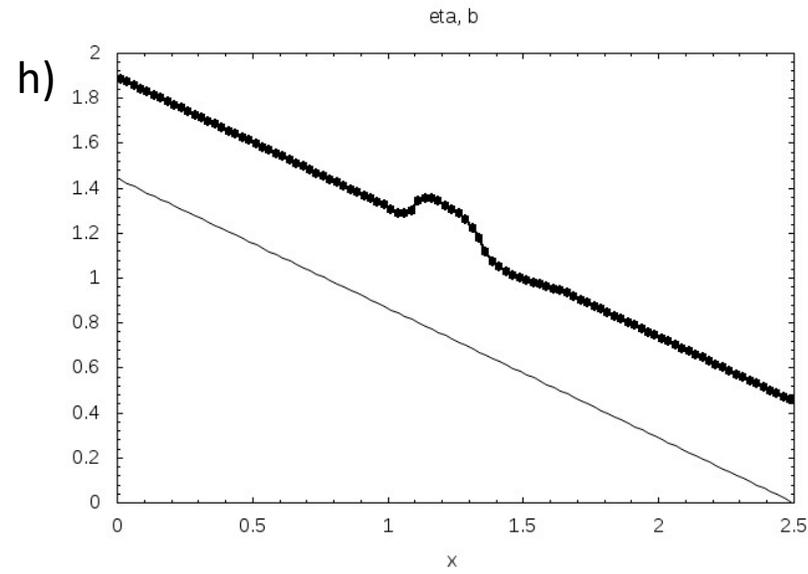 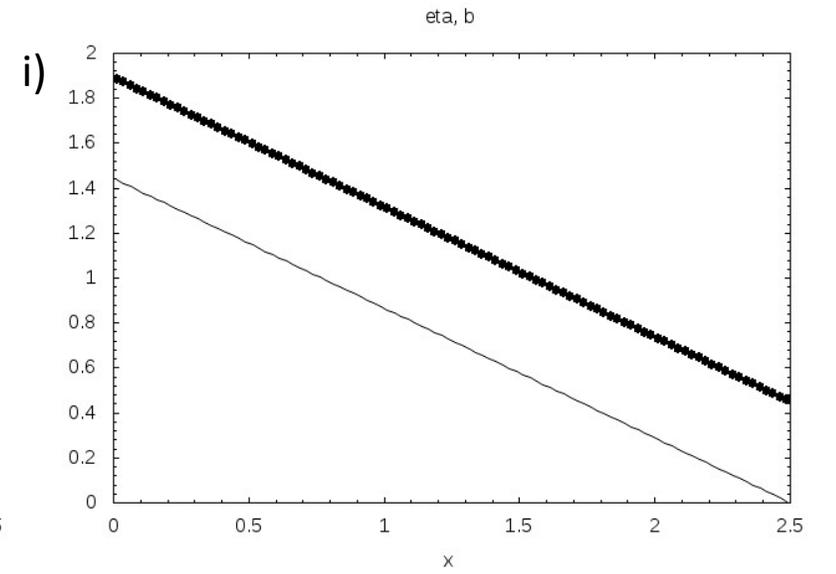

*Fig. 5 shows the free surface and the topography for the shallow water tests for flow over a slanted surface. Figs. 5a, 5b and 5c show the profiles for the test 1 at time $t = 0$, $t = 1$ and $t = 100$, respectively. Figs. 5d, 5e and 5f show the profiles for the test 2 at time $t = 0$, $t = 0.5$ and $t = 100$, respectively. Figs. 5g, 5h and 5i show the profiles for the test 3 at time $t = 0$, $t = 0.05$ and $t = 100$, respectively.*

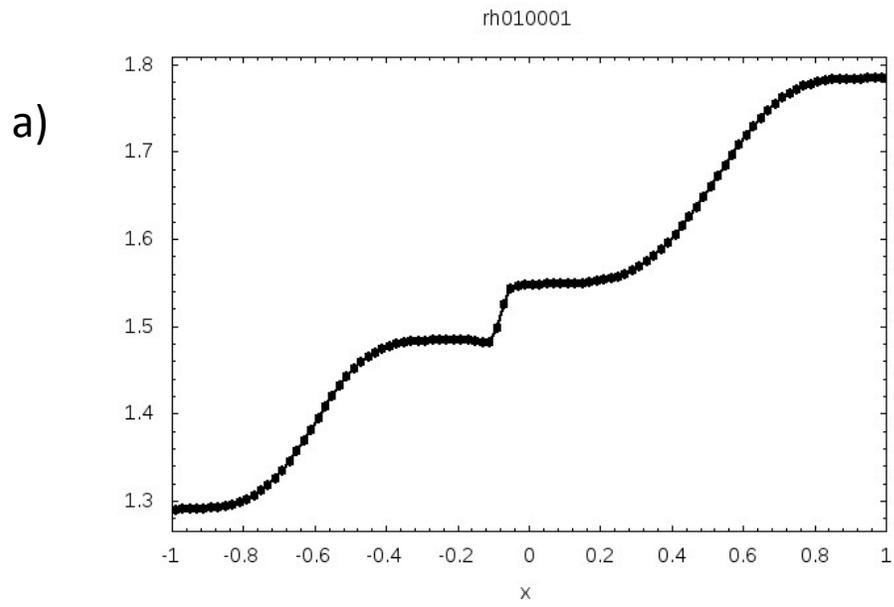
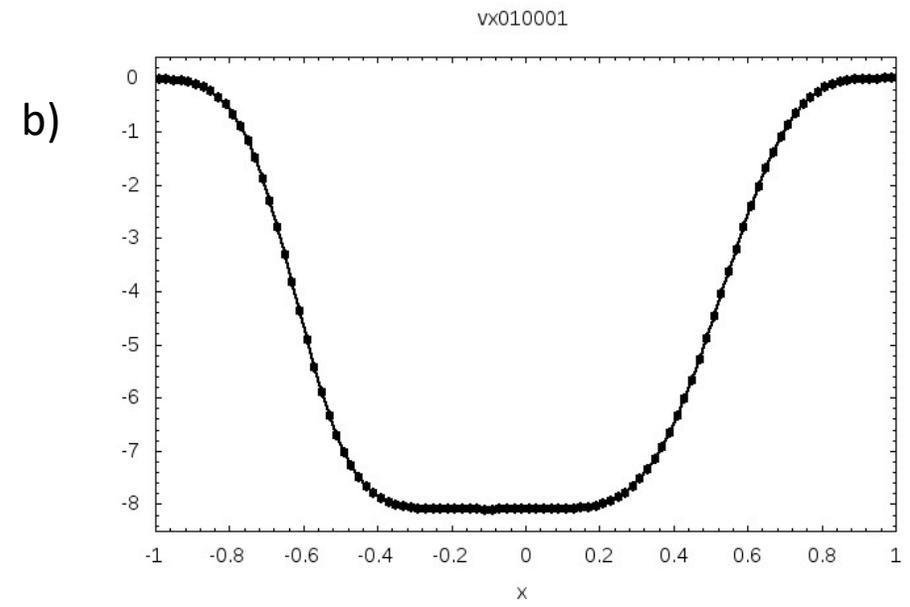
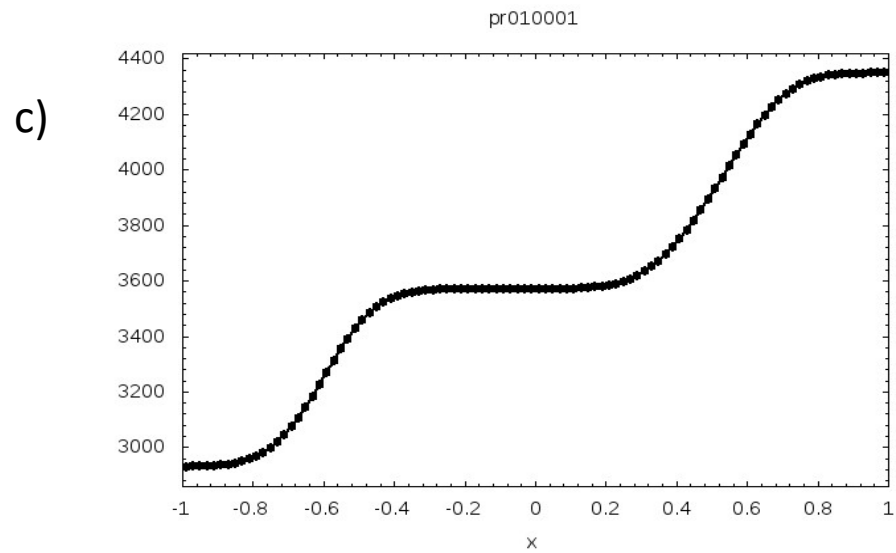
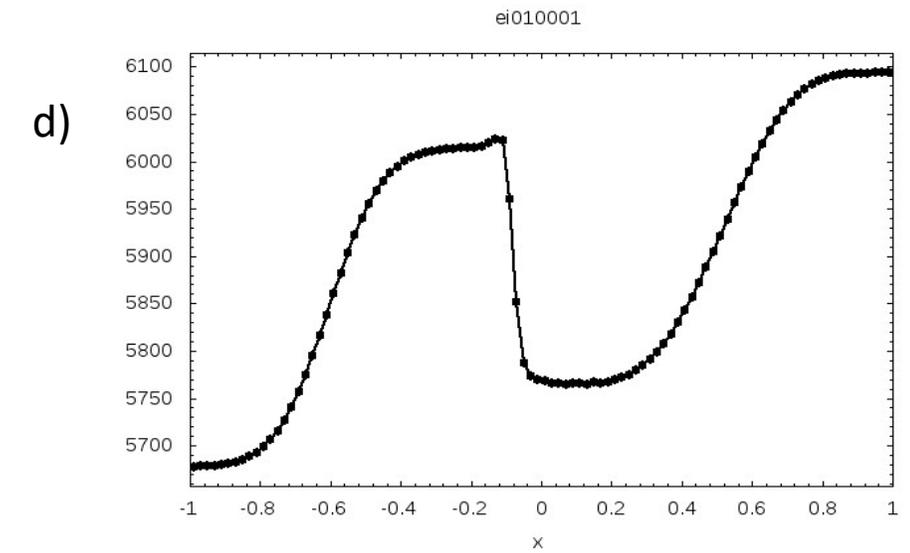

*Fig. 6 shows the results for the compressible Navier-Stokes flow with $\mu = 2Pa/s$ at the output time 0.01. Fig. 6a shows the profile for the density. Fig. 6b shows the profile for the x-component of the velocity. Fig. 6c shows the profile for the density. Fig. 6d shows the internal energy.*

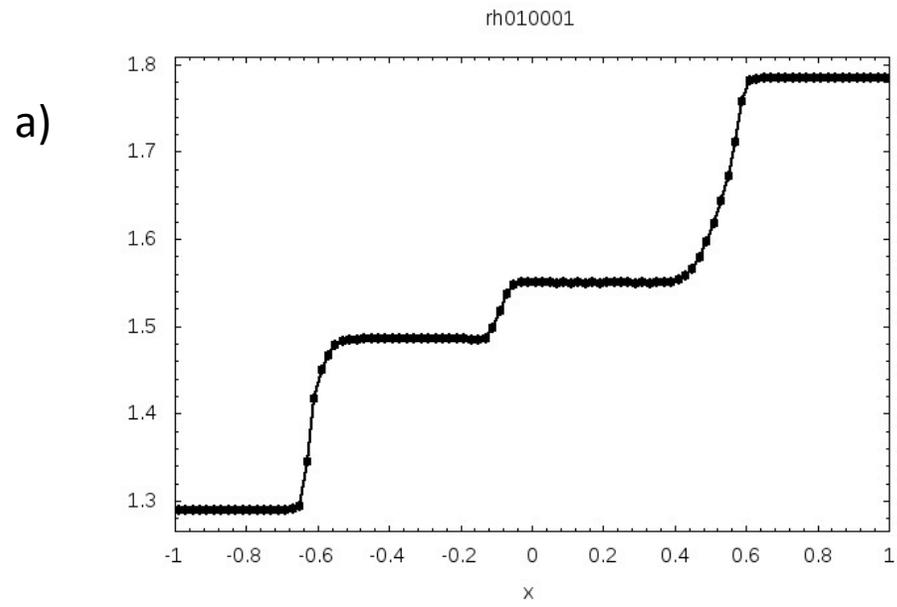
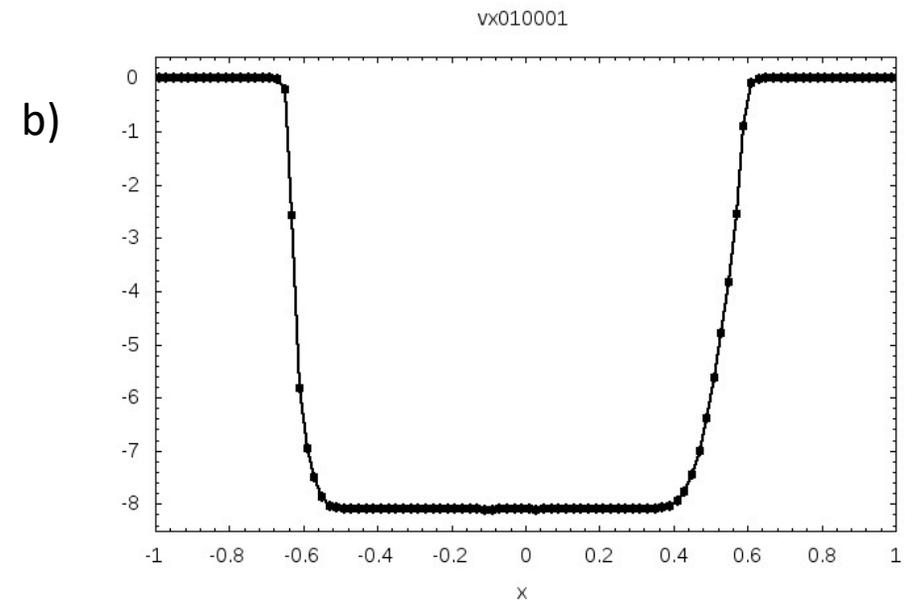
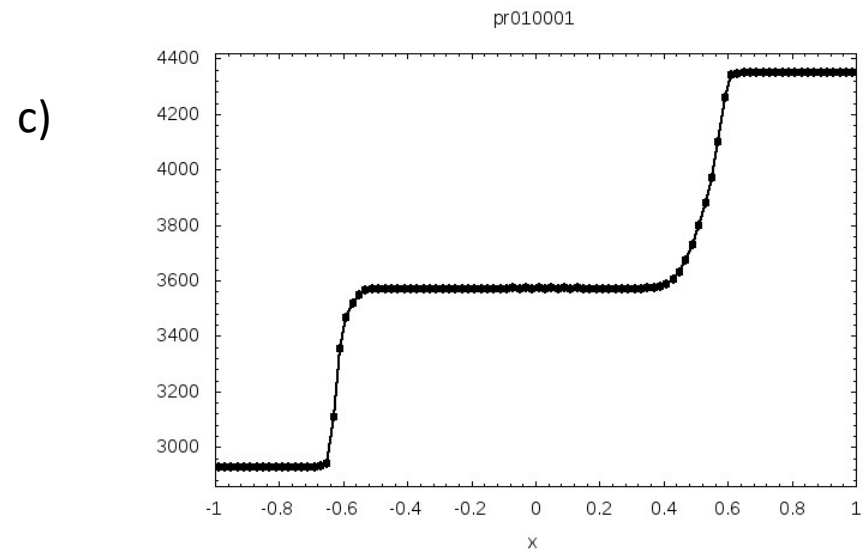
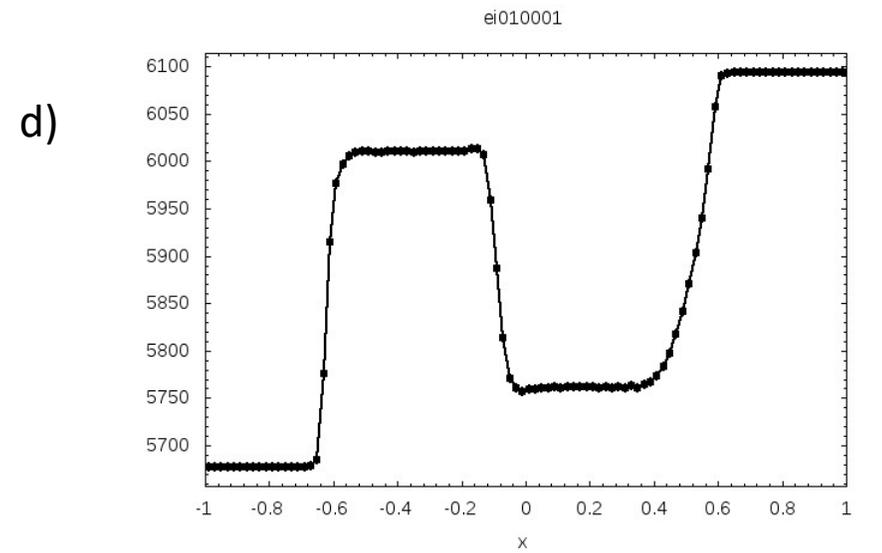

*Fig. 7 shows the results for the compressible Navier-Stokes flow with µ = 0.2Pa/s at the output time 0.01. Fig. 7a shows the profile for the density. Fig. 7b shows the profile for the x-component of the velocity. Fig. 7c shows the profile for the density. Fig. 7d shows the internal energy.*

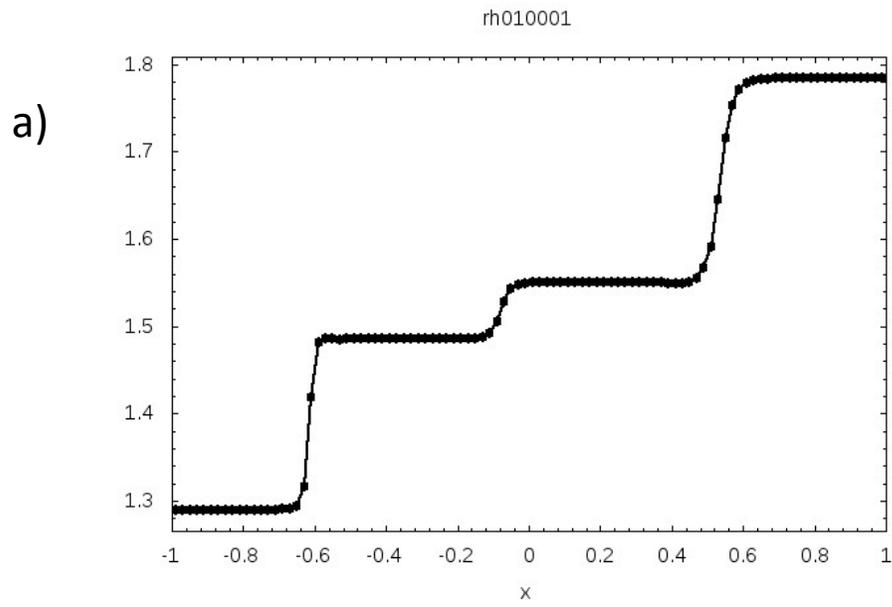 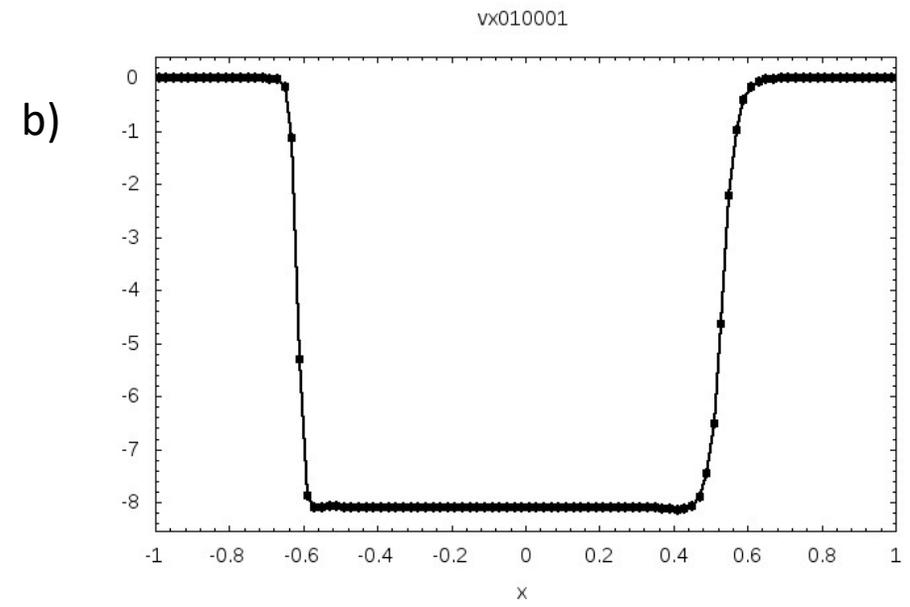 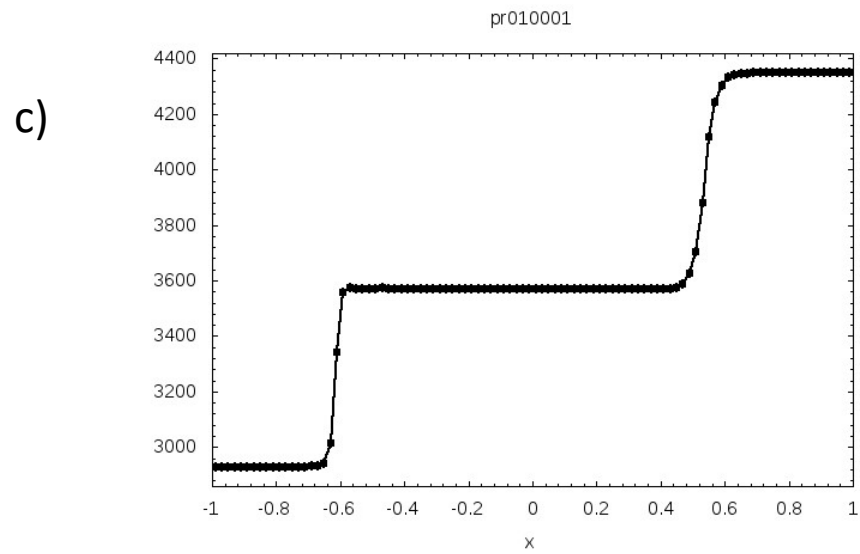 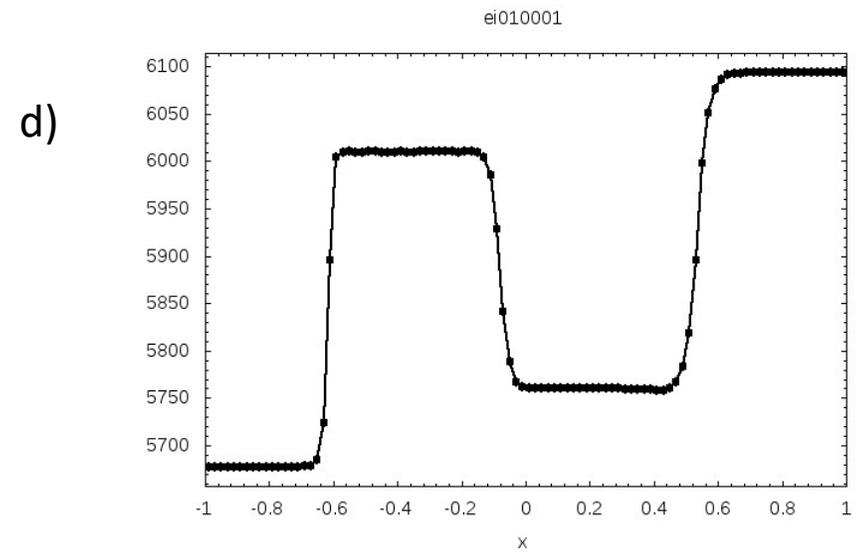

*Fig. 8 shows the results for the compressible Navier-Stokes flow with $\mu = 0.001 Pa/s$ at the output time 0.01. Fig. 8a shows the profile for the density. Fig. 8b shows the profile for the x-component of the velocity. Fig. 8c shows the profile for the density. Fig. 8d shows the internal energy.*